%% A manuscript in PLAIN TeX
\def\author{Kleiman and Piene}
\def\title{Node polynomials for families}
\def\date{\today}
\def\abstract{%
  We continue the development of methods for enumerating nodal curves on
smooth complex surfaces, stressing the range of validity.  We illustrate
the new methods in three important examples.  First, for up to eight
nodes, we confirm G\"ottsche's conjecture about plane curves of low
degree.  Second, we justify Vainsencher's enumeration of irreducible
six-nodal plane curves on a general quintic threefold in four-space.
Third, we supplement Bryan and Leung's enumeration of nodal curves in a
given homology class on an Abelian surface of Picard number one.
 }
%%%%%%%%%%%%%%%%%  \PaperSize is used to center text on page
%\def\PaperSize{letter}
  \def\PaperSize{AFour}         %  center text on page
%%%%%%%%%%%%%%%%%       END of CHOICES
%%%%%%%%%%%%%%%%%       Left-Right
   % to put page on right
% \def\FirstPageOnRight{TRUE}

\def\TheMagstep{\magstep0}
 \def\eqno{\leqno\indent}
%\input b.mac
%%%%%%%%% bol.mac = a Plain TeX macro package by SLK
\let\@=@  %% For e-mail address

 %% For Repeating Similar Acts Defining Macros
\def\GetNext#1 {\def\NextOne{#1}\if\relax\NextOne\let\next=\relax
        \else\let\next=\DoIt \fi \next}
\def\DoIt{\Act\NextOne\GetNext}
\def\ActOn#1{\expandafter\GetNext #1\relax\ }
\def\defcs#1{\expandafter\xdef\csname#1\endcsname}

 %% PAGE LAYOUT
% NB: ams-spec.sty page dimensions = 30pc x 50.5pc
 %Crelle's dimensions
\hsize14.5cm
\vsize22cm
%\parskip=\bigskipamount
 \parindent2pc  %10pt
\abovedisplayskip 4pt plus3pt minus1pt
\belowdisplayskip=\abovedisplayskip
\abovedisplayshortskip 2.5pt plus2pt minus1pt
\belowdisplayshortskip=\abovedisplayskip

\ifx\TheMagstep\undefined \def\TheMagstep{\magstep1} \fi
\def\TRUE{TRUE} % For Boolean tests
\ifx\DoublepageOutput\TRUE \def\TheMagstep{\magstep0} \fi
\mag=\TheMagstep

% CENTER TEXT ON PAGE
        % additional vertical adjustment
\newskip\vadjustskip \vadjustskip=0.5\normalbaselineskip
\def\centertext
 {\hoffset=\pgwidth \advance\hoffset-\hsize
  \advance\hoffset-2truein \divide\hoffset by 2\relax
  \voffset=\pgheight \advance\voffset-\vsize
  \advance\voffset-2truein \divide\voffset by 2\relax
  \advance\voffset\vadjustskip
 }
\newdimen\pgwidth\newdimen\pgheight
\def\letter{letter}\def\AFour{AFour}
\ifx\PaperSize\letter
 \pgwidth=8.5truein \pgheight=11truein
 \message{- Got a paper size of letter.  }\centertext
\fi
\ifx\PaperSize\AFour
 \pgwidth=210truemm \pgheight=297truemm
 \message{- Got a paper size of AFour.  }\centertext
\fi

%% HEADLINE STYLE
\def\today{\ifcase\month\or     % From the TeX book p. 406
 January\or February\or March\or April\or May\or June\or
 July\or August\or September\or October\or November\or December\fi
 \space\number\day, \number\year}
\nopagenumbers
 \newcount\pagenumber \pagenumber=1
 \def\advancepagenumber{\global\advance\pagenumber by 1}
\def\folio{\number\pagenum} % \pagenum is \let below by an \ifx
\headline={%
  \ifnum\pagenum=0\hfill
  \else
   \ifnum\pagenum=1\firstheadline
   \else
     \ifodd\pagenum\oddheadline
     \else\evenheadline\fi
   \fi
  \fi
}
\expandafter\ifx\csname date\endcsname\relax \let\dato=\today
            \else\let\dato=\date\fi
 \message{\dato}
\font\eightrm=cmr8
\let\firstheadline\hfill
\def\oddheadline{% \rlap{\dato}
 \hfil{\eightpoint\it\author, \headtitle}\hfil
 \llap{\rm\folio}}
\def\evenheadline{\rlap{\rm\folio}
  \hfil{\eightpoint\it\author, \headtitle}\hfil
 \llap{\rm\dato}}
\def\headtitle{\title}

%% TWO-COLUMN LANDSCAPE FORMAT
% Modified from the TeX book, p. 257.
 \newdimen\fullhsize \newbox\leftcolumn
 \def\fulline{\hbox to \fullhsize}
\def\doublepageoutput
{\let\lr=L
 \output={\if L\lr
           \global\setbox\leftcolumn=\columnbox \global\let\lr=R%
          \else \doubleformat \global\let\lr=L
          \fi
        \ifnum\outputpenalty>-20000 \else\dosupereject\fi
        }%
 \def\doubleformat{\shipout\vbox{%
     \ifx\PaperSize\AFour
           \fulline{\hfil\box\leftcolumn\hfil\columnbox\hfil}%
     \else
           \fulline{\hfil\hfil\box\leftcolumn\hfil\columnbox\hfil\hfil}%
     \fi             }%
     \advancepageno
}
 \def\columnbox{\vbox
   {\if E\topmark\headline={\hfil}\nopagenumbers\fi
    \makeheadline\pagebody\makefootline\advancepagenumber}%
   }%
\fullhsize=\pgheight \hoffset=-1truein
 \voffset=\pgwidth \advance\voffset-\vsize
  \advance\voffset-2truein \divide\voffset by 2
  \advance\voffset\vadjustskip
 
 %%% to put page on right
\ifx\FirstPageOnRight\TRUE % to put page on right
 \null\vfill\nopagenumbers\eject\pagenum=1\relax
\fi
}
\ifx\DoublepageOutput\TRUE \let\pagenum=\pagenumber\doublepageoutput
 \else \let\pagenum=\pageno \fi

% ADDITIONAL FONTS
% \font\twelvebf=cmbx12          % For title
 \font\seventeenrm=cmr17
 \font\smc=cmcsc10              % For authors' names and section headings
%% EIGHT POINT TYPE for footnotes and references
\catcode`\@=11          % make @ a letter temporarily
\def\eightpoint{\eightpointfonts
 \setbox\strutbox\hbox{\vrule height7\p@ depth2\p@ width\z@}%
 \eightpointparameters\eightpointfamilies
 \normalbaselines\rm
 }
\def\eightpointparameters{%
 \normalbaselineskip9\p@
 \abovedisplayskip4.5\p@ plus2.4\p@ minus1.5\p@
 \belowdisplayskip=\abovedisplayskip
 \abovedisplayshortskip\z@ plus2.4\p@
 \belowdisplayshortskip5.6\p@ plus2.4\p@ minus3.2\p@
 }
\newfam\smcfam
\def\eightpointfonts{%
 \font\eightrm=cmr8 \font\sixrm=cmr6
 \font\eightbf=cmbx8 \font\sixbf=cmbx6
 \font\eightit=cmti8
 \font\eightsmc=cmcsc8
 \font\eighti=cmmi8 \font\sixi=cmmi6
 \font\eightsy=cmsy8 \font\sixsy=cmsy6
 \font\eightsl=cmsl8 \font\eighttt=cmtt8}
\def\eightpointfamilies{%
 \textfont\z@\eightrm \scriptfont\z@\sixrm  \scriptscriptfont\z@\fiverm
 \textfont\@ne\eighti \scriptfont\@ne\sixi  \scriptscriptfont\@ne\fivei
 \textfont\tw@\eightsy \scriptfont\tw@\sixsy \scriptscriptfont\tw@\fivesy
 \textfont\thr@@\tenex \scriptfont\thr@@\tenex\scriptscriptfont\thr@@\tenex
 \textfont\itfam\eightit        \def\it{\fam\itfam\eightit}%
 \textfont\slfam\eightsl        \def\sl{\fam\slfam\eightsl}%
 \textfont\ttfam\eighttt        \def\tt{\fam\ttfam\eighttt}%
 \textfont\smcfam\eightsmc      \def\smc{\fam\smcfam\eightsmc}%
 \textfont\bffam\eightbf \scriptfont\bffam\sixbf
   \scriptscriptfont\bffam\fivebf       \def\bf{\fam\bffam\eightbf}%
 \def\rm{\fam0\eightrm}%
% \tt \ttglue=0.5em plus0.25em minus0.15em
 }

%% Modification of the PLAIN footnote macro for 8pt

\def\footnoterule{\kern-3pt
  \hrule width 2\parindent \kern 2.6pt}
\def\vfootnote#1{\insert\footins\bgroup
 \eightpoint\catcode`\^^M=5\leftskip=0pt\rightskip=\leftskip%% only change
 \interlinepenalty\interfootnotelinepenalty
  \splittopskip\ht\strutbox % top baseline for broken footnotes
  \splitmaxdepth\dp\strutbox \floatingpenalty\@MM
  \leftskip\z@skip \rightskip\z@skip \spaceskip\z@skip \xspaceskip\z@skip
  \textindent{#1}\footstrut\futurelet\next\fo@t}

%%  ``Ties'' with a \thinspace for page numbers
\def\p.{p.\penalty\@M \thinspace}
\def\pp.{pp.\penalty\@M \thinspace}
%% SECTIONING
\newcount\sctno
\def\sctn#1\par
  {\removelastskip\vskip0pt plus4\normalbaselineskip \penalty-250
  \vskip0pt plus-5\normalbaselineskip \bigskip\bigskip\medskip
  \centerline{\bf #1}\nobreak\medskip
}

\def\sct#1 {\sctno=#1\relax\sctn#1. }

%%  STYLE MACROS  %%
%% Redefine \item to give greater indentation than AMSTeX
%   and the roman font within parentheses.
\def\item#1 {\par\indent
 \hangindent\parindent
 \llap{\rm (#1)\enspace}\ignorespaces}
%% item for Crelle
 \def\rbit#1 {\par{\rm (#1)\enspace}\ignorespaces}
%% Define a similar macro without the hanging indentation for assertions
%% and that starts each part with an ordinary \parindent
 \def\inpart#1 {{\rm (#1)\enspace}\ignorespaces}
 \def\part {\par\inpart}

%% ARTICLES
\def\Cs#1){\(\number\sctno.#1)}
\def\part#1 {\par\(#1)\enspace\ignorespaces}

\def\dsc#1 #2.{\medbreak{\bf\Cs#1)} {\it #2.} \ignorespaces}
%% For setting results
\def\proclaim#1 #2 {\medbreak%\goodbreak
  {\bf#1 (\number\sctno.#2).}\quad \bgroup%\activeleftp
\it}
\def\endproclaim{\par\egroup\medbreak}
\def\pf{\endproclaim{\bf Proof.}\quad\ignorespaces}
\def\boxit#1{\vbox{\hrule\hbox{\vrule\kern3pt
	\vbox{\kern3pt#1\kern3pt}\kern3pt\vrule}\hrule}}
\def\qed{\quad\lower1pt\hbox{\boxit\null}}
\def\lem{\proclaim Lemma } 
  \def\thm{\proclaim Theorem }
\def\xpln#1 #2 {\medbreak {\bf #1 (\number\sctno.#2)}.\quad}

\def\rmk{\xpln Remark }
\def\eg{\xpln Example }
%% REFERENCING
        % to introduce the keys in order
 \newcount\refno \refno=0        \def\NoKey{*!*}
 \def\MakeKey{\advance\refno by 1 \expandafter\xdef
  \csname\TheKey\endcsname{{\number\refno}}\NextKey}
 \def\NextKey#1 {\def\TheKey{#1}\ifx\TheKey\NoKey\let\next\relax
  \else\let\next\MakeKey \fi \next}
 \def\RefKeys #1\endRefKeys{\expandafter\NextKey #1 *!* }
 \def\SetRef#1 #2,{\hang\llap
  {[\csname#1\endcsname]\enspace}{\smc #2},}
 \newbox\keybox \setbox\keybox=\hbox{[25]\enspace}
 \newdimen\keyindent \keyindent=\wd\keybox
\def\references{\kern-\medskipamount
  \sctn References\par
  \vskip-\medskipamount
  \bgroup   \frenchspacing   \eightpoint
   \parindent=\keyindent  \parskip=\smallskipamount
   \everypar={\SetRef}\par}
\def\endreferences{\egroup}

%% SERIALS
 \def\serial#1#2{\expandafter\def\csname#1\endcsname ##1 ##2 ##3
        {\unskip\ {\it #2\/} {\bf##1} (##2), ##3}} % \serial{}{}

%% modified \cite code from AMSTeX
\def\UThin{\penalty\@M \thinspace\ignorespaces}
        % unbreakable \thinspace for use after periods
\def\(#1){{\let~=\UThin\rm(#1)}}
\def\relaxnext@{\let\next\relax}
\def\cite#1{\relaxnext@
 \def\nextiii@##1,##2\end@{\unskip\space{\rm[\SetKey{##1}],\let~=\UThin##2}}%
 \in@,{#1}\ifin@\def\next{\nextiii@#1\end@}\else
 \def\next{{\rm[\SetKey{#1}]}}\fi\next}
\newif\ifin@
\def\in@#1#2{\def\in@@##1#1##2##3\in@@
 {\ifx\in@##2\in@false\else\in@true\fi}%
 \in@@#2#1\in@\in@@}
\def\SetKey#1{{\bf\csname#1\endcsname}}

\catcode`\@=12  %at signs are no longer letters, but other

 %%  MATH MACROS
\let\:=\colon \let\ox=\otimes \let\x=\times
 \let\wh=\widehat \let\?=\overline

 \def\onto{\to\mathrel{\mkern-15mu}\to}
 
\def\smashedlongrightarrow{\setbox0=\hbox{$\longrightarrow$}\ht0=1pt\box0}
\def\risom{\buildrel\sim\over{\smashedlongrightarrow}}
 \def\lgto{-\mathrel{\mkern-10mu}\to}
 \def\smashedlgto{\setbox0=\hbox{$\scriptstyle\lgto$}\ht0=1.85pt
        \lower1.25pt\box0}

\def\tsum{\textstyle\sum}
  \def\es{^{\rm es}}

\def\Act#1{\defcs{c#1}{{\fam2#1}}}               % Script
 \ActOn{C F I J K L M N O P Q R }
 \def\Act#1{\defcs{#1}{\mathop{\rm#1}\nolimits}} % Roman operators
 \ActOn{cod sts frs frs$_1$ Hilb rts type wt Supp rk Sat mod ch Pic
  }
\def\Act#1{\defcs{c#1}{\mathop{\it#1}\nolimits}}% Italic operators
 \ActOn{Aut Cok Ext Hom Ker Sym }
\def\Act#1{\defcs{#1}{\hbox{ #1 }}}         % Text for math displays
 \ActOn{and by for where with on }
\def\Act#1{\defcs{I#1}{{\bf#1}}}                % Bold
 \ActOn{A B C D E F G H P Q R S T V X Z m }

%%% End of the macro file

%\input /home/kleiman/art/macro/amscd.sty%%%%%%%%%%%%%%%%%%
%%    7-10-97  10:03p        4,825  amscd.mac
%%%%%%%%%%%%%%%%%%
\let\at=@\catcode`\@=11

 \def\activeat#1{\csname @#1\endcsname}
 \def\def@#1{\expandafter\def\csname @#1\endcsname}
 {\catcode`\@=\active \gdef@{\activeat}}

\let\ssize\scriptstyle
\newdimen\ex@	\ex@.2326ex

 \def\requalfill{\cleaders\hbox{$\mkern-2mu\mathord=\mkern-2mu$}\hfill
  \mkern-6mu\mathord=$}
 \def\eqfill{$\m@th\mathord=\mkern-6mu\requalfill}
 \def\deffill{\hbox{$:=$}$\m@th\mkern-6mu\requalfill}
 \def\fiberbox{\hbox{$\vcenter{\hrule\hbox{\vrule\kern1ex
     \vbox{\kern1.2ex}\vrule}\hrule}$}}

 \font\arrfont=line10
 \def\Swarrow{\vcenter{\hbox{$\swarrow$\kern-.26ex
    \raise1.5ex\hbox{\arrfont\char'000}}}}

 \newdimen\arrwd
 \newdimen\minCDarrwd \minCDarrwd=2.5pc
 	
 \def\findarrwd#1#2#3{\arrwd=#3%
  \setbox\z@\hbox{$\ssize\;{#1}\;\;$}%
  \setbox\@ne\hbox{$\ssize\;{#2}\;\;$}%
  \ifdim\wd\z@>\arrwd \arrwd=\wd\z@\fi
  \ifdim\wd\@ne>\arrwd \arrwd=\wd\@ne\fi}
 \newdimen\arrowsp\arrowsp=0.375em
 \def\findCDarrwd#1#2{\findarrwd{#1}{#2}{\minCDarrwd}
    \advance\arrwd by 2\arrowsp}
 \newdimen\minarrwd 
 \setbox\z@\hbox{$\longrightarrow$} \minarrwd=\wd\z@

 \def\harrow#1#2#3#4{{\minarrwd=#1\minarrwd
   \findarrwd{#2}{#3}{\minarrwd}\kern\arrowsp
    \mathrel{\mathop{\hbox to\arrwd{#4}}\limits^{#2}_{#3}}\kern\arrowsp}}

%%%%% -SLK-00Jun25> Copied from nlcd
%***arrow heights:
\newdimen\Bigheight
\setbox0=\hbox{$\Big\downarrow$}\Bigheight=\ht0
\newdimen\arrowheight
\arrowheight=\Bigheight\advance\arrowheight by \dp0
\def\CDBig#1{{\hbox{$\left#1\vbox to\Bigheight{}\right.\n@space$}}}

 \def@]#1>#2>#3>{\harrow{#1}{#2}{#3}\rightarrowfill}
 \def@>#1>#2>{\harrow1{#1}{#2}\rightarrowfill}
 \def@<#1<#2<{\harrow1{#1}{#2}\leftarrowfill}
 \def@={\harrow1{}{}\eqfill}
 \def@:#1={\harrow1{}{}\deffill}
 \def@ N#1N#2N{\vCDarrow{#1}{#2}\UpDownarrow}
 \def\UpDownarrow{\uparrow\,\Big\downarrow}

%%  -SLK-94JUNE14=> hooked horizontal arrows
\def\hookrightarrowfill{\hbox{$\lhook\joinrel$}\rightarrowfill}
\def@{hk>}#1>#2>{\harrow1{#1}{#2}\hookrightarrowfill}
\def@ c>#1>#2>{\harrow1{#1}{#2}\hookrightarrowfill}
\def\hookleftarrowfill{\leftarrowfill\hbox{$\joinrel\rhook$}}
\def@{hk<}#1<#2<{\harrow1{#1}{#2}\hookleftarrowfill}

 \def@.{\ifodd\row\relax\harrow1{}{}\hfill
   \else\vCDarrow{}{}.\fi}
 \def@|{\vCDarrow{}{}\Vert}
 \def@ V#1V#2V{\vCDarrow{#1}{#2}\downarrow}
 \def@ A#1A#2A{\vCDarrow{#1}{#2}\uparrow}
 \def@(#1){\arrwd=\csname col\the\col\endcsname\relax
   \hbox to 0pt{\hbox to \arrwd{\hss$\vcenter{\hbox{$#1$}}$\hss}\hss}}

 \def\squash#1{\setbox\z@=\hbox{$#1$}\finsm@@sh}
\def\finsm@@sh{\ifnum\row>1\ht\z@\z@\fi \dp\z@\z@ \box\z@}

 \newcount\row \newcount\col \newcount\numcol \newcount\arrspan
 \newdimen\vrtxhalfwd  \newbox\tempbox

 \def\innernewdimen{\alloc@1\dimen\dimendef\insc@unt}
 \def\measureinit{\col=1\vrtxhalfwd=0pt\arrspan=1\arrwd=0pt
   \setbox\tempbox=\hbox\bgroup$}
 \def\setinit{\col=1\hbox\bgroup$\ifodd\row
   \kern\csname col1\endcsname
   \kern-\csname row\the\row col1\endcsname\fi}
 \def\findvrtxhalfsum{$\egroup
  \expandafter\innernewdimen\csname row\the\row col\the\col\endcsname
  \global\csname row\the\row col\the\col\endcsname=\vrtxhalfwd
  \vrtxhalfwd=0.5\wd\tempbox
%  \global\expandafter\advance\csname row\the\row col\the\col\endcsname
  \global\advance\csname row\the\row col\the\col\endcsname by \vrtxhalfwd
  \advance\arrwd by \csname row\the\row col\the\col\endcsname
  \divide\arrwd by \arrspan
  \loop\ifnum\col>\numcol \numcol=\col%
     \expandafter\innernewdimen \csname col\the\col\endcsname
     \global\csname col\the\col\endcsname=\arrwd
   \else \ifdim\arrwd >\csname col\the\col\endcsname
      \global\csname col\the\col\endcsname=\arrwd\fi\fi
   \advance\arrspan by -1 %
   \ifnum\arrspan>0 \repeat}
 \def\setCDarrow#1#2#3#4{\advance\col by 1 \arrspan=#1
    \arrwd= -\csname row\the\row col\the\col\endcsname\relax
    \loop\advance\arrwd by \csname col\the\col\endcsname
     \ifnum\arrspan>1 \advance\col by 1 \advance\arrspan by -1%
     \repeat
    \squash{\mathop{
     \hbox to\arrwd{\kern\arrowsp#4\kern\arrowsp}}\limits^{#2}_{#3}}}
 \def\measureCDarrow#1#2#3#4{\findvrtxhalfsum\advance\col by 1%
   \arrspan=#1\findCDarrwd{#2}{#3}%
    \setbox\tempbox=\hbox\bgroup$}
 \def\vCDarrow#1#2#3{\kern\csname col\the\col\endcsname
    \hbox to 0pt{\hss$\vcenter{\llap{$\ssize#1$}}%
     \Big#3\vcenter{\rlap{$\ssize#2$}}$\hss}\advance\col by 1}

 \def\setCD{\def\harrow{\setCDarrow}%
  \def\\{$\egroup\advance\row by 1\setinit}
  \m@th\lineskip3\ex@\lineskiplimit3\ex@ \row=1\setinit}
 \def\endsetCD{$\egroup}
 \def\drop#1\\{\findvrtxhalfsum\advance\row by 2 \measureinit}
 \def\measure{\bgroup
  \def\harrow{\measureCDarrow}%
   \def\\##1{\ifx##1\endmeasure\endmeasure\else\expandafter\drop\fi}
 \row=1\numcol=0\measureinit}
 \def\endmeasure{\findvrtxhalfsum\egroup}

\newbox\CDbox \newdimen\sdim

 \newcount\savedcount
 \def\CD#1\endCD{\savedcount=\count11%
   \measure#1\endmeasure
%   \vcenter{\setCD#1\endsetCD}%
%   \vcenter{\setCD#1\endsetCD\kern\smallskipamount}%
   \vcenter{\setCD#1\endsetCD\kern\medskipamount}%
%   \vbox{\setCD#1\endsetCD}%
%   \vcenter spread\sdim {\setCD#1\endsetCD}%
%   \global\setbox\CDbox\vbox spread\sdim {\setCD#1\endsetCD}%
%   \copy\CDbox
   \global\count11=\savedcount}

 \catcode`\@=\active

%%%%%%%%%%%%%%%%%       REFERENCE KEYS
   \RefKeys
 BS88 B--L99 BL98 CH98 C66 C86 Co70 Ch97 CK99 D99 Gtt98 GL96 GLS97 HP95
JK96 TdJ00 Ka86 Ka92 Kl74 KP99 K--P K01 Los98 Mat91 M70 NV97 Ran89 R67
R75 S65 S79 St48 Ta82 Te73 V95 V97 Wahl74 Wall84 Z82
   \endRefKeys
%%%%%%%%%%%%%%%%%       TOPMATTER

\long\def\TOPMATTER
{\null\vskip3.5true cm
%   \obeylines
   {\seventeenrm
	 \centerline{Node polynomials for families:}\smallskip
	\centerline{results and examples}
   }  \bigskip
 %% Subject classification and acknowledgements
   \footnote{}{\noindent %
   MSC-class: 14N10 (Primary); 14C20, 14H40, 14K05 (Secondary).}
 %% Authors' names, addresses and support
\centerline{By {\it Steven Kleiman\/} and {\it Ragni Piene\/}}
  \bigskip
  \centerline{\vrule width 3 true cm height 0.4pt depth 0pt}
   \medskip
   %% Abstract
 {%\parindent=1.5\parindent \narrower  \noindent
  % \eightpoint{\smc Abstract.}\enspace \ignorespaces\abstract \par}
  {\bf Abstract.}\quad\abstract\par}
 %\bigskip
} %%\end of Topmatter

   %% Body
%  \predisplaypenalty=500 \postdisplaypenalty=-500

 \TOPMATTER
 \kern-\baselineskip
\sct1 Introduction

This paper is the second in a series devoted to the enumeration of nodal
curves on smooth complex surfaces.  The first paper \cite{KP99} focuses
on curves in a ``suitably'' ample linear system on a fixed ambient
surface.  This second paper treats more general systems and variable
surfaces, and places greater stress on the range of validity.  Here we
develop some general methods, and use them in three important examples:
curves of low degree in the plane, plane curves on a threefold in
four-space, and homologous curves on an Abelian surface.

Nodal plane curves were enumerated, for up to three nodes, in the third
quarter of the nineteenth century, and the general problem has recently
been revived; for some pertinent history, see Remark~(3.7).  In
particular, G\"ottsche conjectured in \cite{Gtt98, Conj.~4.1, p.~530},
that, for each $r$, if $N_r(m)$ denotes the number of curves of degree
$m$ with $r$ nodes through $m(m+3)/2-r$ general points, then $N_r(m)$ is
given by a certain ``node'' polynomial of degree $2r$ in $m$ for $m\ge
r/2+1$, which is just the range of $m$ where the locus of nonreduced
curves is too small to interfere.  Our first main result, Theorem~(3.1),
confirms G\"ottsche's conjecture for $r\le8$.

Theorem~(3.1) can be derived from Theorem~(1.3) of \cite{KP99} and the
recursive enumerative formula of Caporaso and Harris \cite{CH98,
p.~353}; see the end of Remark~(3.7).  However, we proceed differently
for three reasons.  First, our approach is self-con\-tained.
 Second, our approach may eventually lead to a confirmation of
G\"ott\-sche's conjecture for all $r$, whereas the alternative approach
requires evaluating Caporaso and Harris's formula at least once for each
$r$, an absurd project. Third, our lemmas are also needed to prove our
second main result, Theorem~(4.1).

Theorem~(4.1) enumerates, for $m\ge4$, the 6-nodal plane curves of
degree $m$ on a general threefold of degree $m$ in 4-space, or what is
the same, its 6-tangent 2-planes.  This enumeration provides a nice
example of the use of our machinery in the case of a nontrivial family
of ambient surfaces.  The family consists of all the planes in 4-space,
parameterized by the Grassmann variety; so each surface is the same, but
the family is nonconstant.  The curves are those cut out on the planes
by the threefold.  The number of curves is given by a certain ``node''
polynomial of degree 18 in $m$.

This enumeration was originally done by Vainsencher \cite{V95}.  Indeed, his
paper inspired this one and its companions \cite{KP99} and
\cite{K--P}; our work just refines and extends his.  In this paper,
notably, we develop some new ways of extending the range of validity of
the enumerations.  For example, for plane curves on threefolds,
Vainsencher's Propositions 3.5 and 4.1 imply only that there exists an
undetermined integer $m_0$ such that the enumeration is valid for $m\ge
m_0$, whereas we handle $m\ge4$.

The case $m=5$ is particularly important because of Clemens' conjecture
and mirror symmetry.  Clemens' conjecture \cite{C86, p.~639}, asserts
notably that, on the general quintic threefold, there are only finitely
many rational curves of each degree, and all are smooth.  Their number
was predicted in 1991 in a dramatic application of mirror symmetry, its
first application to enumerative geometry.  This enumeration is
revisited several times in Cox and Katz's lovely text \cite{CK99}.

The irreducible 6-nodal plane quintics are rational, but singular!  So
this part of Clemens' original conjecture is false, and is not made part
of the conjecture's modern formulation \cite{CK99, p.~202}.
Furthermore, mirror symmetry includes these 6-nodal curves in its count.
However, Pandharipande \cite{CK99, (7.54), p.~206}, found something
worse: each 6-nodal curve has six previously unconsidered double covers.
So, in degree 10, mirror symmetry simply produced the wrong number.  It
cannot be the number of all rational curves, smooth and singular!  It is
too large by six times the number of irreducible 6-nodal curves.

The irreducible 6-nodal curves too were originally enumerated by
Vainsencher in \cite{V95, pp.~513--14}, and we recover his number in our
third main result, Theorem~(4.3), pursuing his approach.  Namely, we use
Theorem~(4.1) to obtain the number of all 6-nodal curves, and from it,
we subtract the number of reducible ones.  However, here again, we
advance Vainsencher's work by paying careful attention to the validity
of the numbers involved.

Our fourth and last main result, Theorem~(5.2), enumerates the
irreducible curves having $r$ nodes and lying in a given homology class
$\gamma$ on an Abelian surface $A$ with Picard number 1.  Say $\gamma$
has self-intersection number $d$, and is $m$ times the positive
primitive class.  Set $g:=d/2-r+1$, and let $N_{g,r}$ be the number of
curves through $g$ general points.  Theorem~(5.2) asserts that, if $r\le
8$, then $N_{g,r}$ is given by a certain polynomial of degree $r+1$ in
$g$ for $g>g_0$ where $g_0$ is a certain number depending on $m$ and
$r$, but not on $A$.  The nine polynomials are listed in Table~(5.1).

The first theorem of this sort was proved by Bryan and Leung \cite{BL98,
Thm.~1.1, p.~312}, using symplectic methods.  Their theorem is valid for
any $r$ and $g$, provided $A$ is generic in the following sense: given
the underlying topological space, the complex structure of $A$ is
generic among those for which the given class $\gamma$ is algebraic.  It
follows (as stated in the proof of \cite{BL98, Lem.~5}) that $A$ has
Picard number 1 and that $\gamma$ is primitive, that is, $m=1$.  By
contrast, we fix $A$, not $\gamma$; moreover, our methods are
algebraic-geometric and rather different.  Thus our work supplements
theirs.

Bryan and Leung expressed the $N_{g,r}$ essentially as follows:
 $$\sum_{r\ge 0}(N_{g,r}/g)q^r
	= \biggl(\sum_{k\ge 1}k\sigma_1(k)q^{k-1}\biggr)^{g-1}
	\where \sigma_1(k):=\sum_{d|k}d.$$
 Say the logarithms of the left and right sides are $\sum_{r\ge1}
a_rq^r/r!$ and $(g-1)\sum _{r\ge1}b_rq^r/r!$.  Then $a_r=(g-1)b_r$ for
$r\ge1$.  Moreover, $b_1$, \dots, $b_8$ are these integers:
 $$6,\ -12,\ 168,\ -2448,\ 46944,\ -1071360,\ 29064960,\ -921110400.$$
 Furthermore, there is a weighted homogeneous polynomial $P_r$ of degree
$r$ such that
	$$N_{g,r}=gP_r(a_1,\dots,a_r)/r!\,.$$
   The $P_r$ are defined by the formal identity (2.2), and are known as
the {\it Bell\/} polynomials.  They appear in all our enumerations,
although in the case at hand they enter our work somewhat differently.

Closely related is the enumeration of the $r$-nodal curves lying in a
given linear-equivalence subclass and passing through $g-2$ general
points.  Various cases have been discussed by various authors, and their
work is surveyed in Remark~(5.4).  In particular, G\"ottsche conjectured
a generating function similar to the one above for $N_{g,r}$, and Bryan
and Leung confirmed it when $A$ is generic in the above sense.
Supplementing their work, we can modify the proof of our Theorem~(5.1)
to confirm G\"ottsche's conjecture when $r\le8$ and $m>(3r+5)/2$.

All our enumerations are carried out on the basis of Theorem~(2.5).  Its
statement is implicit in Section~4 of \cite{KP99}; its proof is outlined
there, and is completed in \cite{K--P}.  Here, notably we refine our
treatment of the key cycles.  In \cite{KP99}, they are placed under
unnecessarily stringent genericity hypotheses, which conceivably are not
satisfied in the present circumstances.  So we must adopt a more liberal
definition of these cycles, and develop suitable conditions that imply
the cycles have the right support and are reduced.

More precisely, Theorem~(2.5) concerns a smooth, projective family of
surfaces, $\pi:F\to Y$, where $Y$ is equidimensional and
Cohen--Macaulay.  In Section~4 of \cite{KP99}, mistakenly, $Y$ is not
assumed to be Cohen--Macaulay; on the other hand, unnecessarily, $Y$ is
assumed to be reduced, and the surfaces $\pi^{-1}(y)$, to be
irreducible.

Let $D\subset F$ be a $Y$-flat closed family of curves.  Denote its
rational equivalence class by $v$, and the Chern class
$c_i(\Omega^1_{F/Y})$ by $w_i$.  Partition $Y$ into locally closed
subsets: one $Y(\infty)$ where the fibers $D_y$ have a multiple
component, and each other where the $D_y$ have a given equisingularity
type.  Given $r$, assume that, if nonempty, $Y(\infty)$ has codimension
$r+1$ and that each remaining nonempty subset has codimension at least
$\min(r+1,\,c)$ where $c$ is its expected codimension.

Consider the set of $y\in Y$ where the $D_y$ are $r$-nodal.
Theorem~(2.5) asserts that this set is either empty or exactly of
codimension $r\/$; either way, its closure is the support of a natural
nonnegative cycle $U(r)$.  Furthermore, if $r\le8$, then the class
$[U(r)]$ is equal to $P_r(a_1,\ldots,a_r)/r!$ where $P_r$ is the Bell
polynomial, where $a_q:=\pi_*b_q$, and where $b_q$ is a certain
polynomial in the classes $v,\,w_1,\,w_2$.

In order to apply Theorem~(2.5), we must check that the relevant subsets
of $Y$ have appropriate codimensions.  To do so, we modify several
arguments in \cite{KP99}, and thereby obtain better results.  In the
case of an ambient Abelian surface, basically we replace the Gotzmann
regularity theorem and Bertini's theorem by the Beltrametti--Sommese
$k$-very ampleness theorem.  In the case of curves in the plane, we
take a different tack: we work directly on $Y$ using some of Greuel
and Lossen's results about equisingular families of curves.  Finally, in
the case of 6-nodal plane curves on a threefold, we derive what we
need from our work with curves in the plane.

In each case, therefore, Theorem~(2.5) provides us with an enumerating
cycle $U(r)$ and an effective expression for its class $[U(r)]$.  To
complete the enumeration, we must show that $U(r)$ is reduced so that we
know that each $r$-nodal curve is counted with multiplicity 1.  We do so
by carrying a bit further our analysis of the relevant subsets of $Y$.
Finally, we need to work out the cycles $b_q$, $a_q$, and $P_r$.  This
work is done in Section~4 of \cite{KP99} for any linear system on any
fixed surface, and so it applies in particular to the case of curves in
the plane.  In the remaining two cases, the details are explained, but
the more mechanical calculations are omitted.

In short, in Section~2, we state the general enumeration theorem,
Theorem~(2.5), and explain its ingredients: the Bell polynomials, the
polynomials giving the $b_q$ in terms of $v,\,w_1,\,w_2$, the key
subsets of $Y$, and the enumerating cycle $U(r)$.  In Sections~3, 4, and
5, we work out in detail the three examples: the plane, a threefold
in four-space, and an Abelian surface.

\sct2 The general theorem

In this section, we discuss the general enumeration theorem, Theorem
(2.5), that we use in the following sections.  In Remark~(2.7), we
conjecture a possible generalization.

Let $Y$ be an equidimensional Cohen--Macaulay scheme of finite type over
the complex numbers.  Let $\pi\: F\to Y$ be a smooth projective family
of surfaces, and $D$ a relative effective divisor on $F/Y$.  Fix
$r\ge0$, and consider the points $y\in Y$ parameterizing the curves
$D_y$ with precisely $r$-nodes.

Theorem (2.5) says that these $y$ are enumerated by a cycle $U(r)$, and
that if $r\le8$, then the rational equivalence class $[U(r)]$ is given
by a universal polynomial in the classes $y(a,b,c)$ that are defined as
follows:
        $$y(a,b,c):=\pi_*v^aw_1^bw_2^c\where v:=[D]
         \and w_i:=c_i(\Omega^1_{F/Y}).\eqno\Cs1)$$
 The only hypotheses are that certain key subsets of $Y$ have
appropriate dimensions.

The universal polynomial has a special shape, which makes it much easier
to find and evaluate.  Namely, define  auxiliary polynomials
$P_i(a_1,\ldots,a_i)$ via this formal identity in $t$:
       $$\textstyle \sum_{i\ge 0}P_it^i/i!
          := \exp\bigl(\sum_{j\ge 1}a_j t^j/j!\bigr).\eqno\Cs2)$$
 For example, $P_0=1$, and $P_1=a_1$, and $P_2=a_1^2+a_2$, and
$P_3=a_1^3+3a_1a_2+a_3$.  If we assign $a_j$ weight $j$, then
$P_i(a_1,\ldots,a_i)$ is weighted homogeneous of degree $i$.  These
polynomials are known as the (complete) {\it Bell\/} polynomials, and
have been studied by a number of authors; see Comtet's book \cite{Co70,
\pp.144--48}.

The universal polynomial can be obtained from $P_r(a_1,\ldots,a_r)/r!$
by replacing each $a_q$ by a certain linear combination of the
$y(a,b,c)$ with $a+b+2c=q+2$.  Equivalently, we can set $a_q:=\pi_*b_q$
where $b_q$ is a certain weighted homogeneous polynomial of degree $q+2$
in $v,\,w_1,\,w_2$ if we assign $v$ and $w_1$ weight 1 and $w_2$ weight
2.  The $b_q$ are given by a simple algorithm; it was stated informally in
Section~4 of \cite{KP99}, and is stated in pseudo-code in Algorithm
(2.3) below.

\goodbreak\midinsert
% \medskip
 \centerline{\bf Algorithm (2.3)}
\centerline{\bf Pseudocode for the  $b_q(v,w_1,w_2)$}
\medskip\hrule\medskip%\smallskip
{\obeylines\parindent=0pt\parskip=0pt
INPUT: indeterminates  $v,\,w_1,\,w_2$.
OUTPUT: polynomials $b_q(v,w_1,w_2)$ for  $q=1,\dots,8$.
\smallskip
FUNCTION:  $Q(i, R)$.\smallskip
{\parindent=2em
INPUT: an integer $i$ and a polynomial  $R(v,w_1,w_2)$.
LOCAL: an indeterminate  $e$.
\smallskip
$R':=R(v-ie,\,w_1+e,\,w_2-e^2)$.
  $R'':=$ the remainder in  $e$ of $R'$ on division by $(e^3+w_1e^2+w_2e)$.
\smallskip
RETURN: $Q(i, R):=-$Coeff$(R'',\,e^2)$.}\smallskip
%END.
\smallskip
 $x_2:=v^3+v^2w_1+vw_2$.
{\openup1\jot
FOR  $s$  FROM 0 TO 2 DO
\qquad $b_{s+1}:=P_s(Q(2,b_1),\ldots,Q(2,b_s))x_2$.
\smallskip %OD.
 $x_3:=v^6+4v^5w_1+5v^4(w_1^2+w_2)+v^3(2w_1^3+11w_1w_2) %
            +v^2(6w_1^2w_2+4w_2^2)+4vw_1w_2^2$.
FOR  $s$  FROM 3 TO 6 DO
  \qquad $b_{s+1}:=P_s(Q(2,b_1),\ldots,Q(2,b_s))x_2$
  \hskip6em$-s(s-1)(s-2)P_{s-3}(Q(3,b_1),\ldots,Q(3,b_{s-3}))x_3$.
\smallskip %OD.
  $x_4:=v^{10}+10v^9w_1+v^8(40w_1^2+15w_2)+v^7(82w_1^3+111w_1w_2)$
  \hskip3em $+v^6(91w_1^4+315w_1^2w_2%
           +63w_2^2)+v^5(52w_1^5+29w_1^3w_2+324w_1w_2^2)$
  \hskip3em $+v^4(12w_1^6+282w_1^4w_2+593w_1^2w_2^2+85w_2^3) %
         +v^3(72w_1^5w_2+464w_1^3w_2^2+259w_1w_2^3)$
  \hskip3em $+v^2(132w_1^4w_2^2+246w_1^2w_2^3+36w_2^4)%
            +v(72w_1^3w_2^3+36w_1w_2^4)$.
 $b_8:= P_7(Q(2,b_1),\ldots,Q(2,b_7))x_2 %
        -7\cdot 6\cdot 5\, P_4(Q(3,b_1),\ldots,Q(3,b_{4}))x_3%
        +3281\cdot 7!\, x_4$.
\par}}
\medskip%\smallskip
\hrule
%\vskip-\normalbaselineskip
\endinsert

Hypothesis (i) of Theorem (2.5) concerns the set $Y(\infty)$ of $y\in Y$
such that the curve $D_y$ has a multiple component, or equivalently, is
nonreduced.  Now, let $\ID$ be a ``minimal Enriques diagram'' as defined
in Section~2 of \cite{KP99}.  Hypothesis (ii) concerns the (locally
closed) set $Y(\ID)$ of $y\in (Y-Y(\infty))$ such that $D_y$ has $\ID$
as its associated diagram.

Briefly put, these $\ID$ are abstract combinatorial structures that
represent the equisingularity types of reduced curves on smooth
surfaces.  The $\ID$ associated to a curve $C$ is made from its directed
resolution graph $\Gamma$.  Weight $\Gamma$ with the multiplicities of
the strict transforms of $C$, and equip $\Gamma$ with the binary
relation of ``proximity''; by definition, one infinitely near point of
$C$ is {\it proximate\/} to a second if the first lies on the strict
transform of the exceptional divisor of the blowup centered at the
second.  By the theorem of embedded resolution, almost all infinitely
near points have multiplicity 1, and are proximate solely to their
immediate predecessors.  Form all the infinite unbroken successions of
these points, and consider the corresponding vertices in the weighted
and equipped $\Gamma$; remove these vertices to get $\ID$.

{}From $\ID$, we can, in principle, determine all the numerical invariants
of the equisingularity class of $C$.  Six such invariants were studied
in Sections~2 and 3 of \cite{KP99}, and they will be used here; so we
recall them now.  Each is given by a formula in these basic numbers:
  $$\eqalign{
  m_V&:=\hbox{the multiplicity, or weight, of the  vertex $V\in\ID$,}\cr
  \frs(\ID)&:=\hbox{the number of free vertices in $\ID$,}\cr
  \rts(\ID)&:=\hbox{the number of roots in $\ID$,}\cr}$$
 where a {\it root\/}  is an initial vertex and a {\it free vertex\/} is one
that is not proximate to a remote predecessor (so a root is free).  Each
remaining vertex is proximate to two vertices, and is said to be a {\it
satellite\/} of the more distant of the two.

 The six numerical invariants are  the following:
  $$\eqalign{
        \dim(\ID)&:=\rts(\ID)+\frs(\ID),\cr
        \deg(\ID)&:=\tsum_{V\in\ID}{m_V+1\choose2},\cr
        \cod(\ID)&:=\deg(\ID)-\dim(\ID),\cr}\quad
 \eqalign{\delta(\ID)&:=\tsum_{V\in\ID}{m_V\choose2},\cr
        r(\ID)&:=\tsum_V \bigl(m_V-\tsum_{W\succ V}m_W\bigr),\cr
        \mu(\ID)&:=2\delta(\ID)-r(\ID)+\rts(\ID),\cr}$$
 where $W\succ V$ means that $W$ is proximate to $V$.  The numbers in
the right column are, respectively, equal to the $\delta$-invariant, the
number of branches, and the Milnor number of $C$.  The numbers in the
left column have geometric meanings, which were discussed in Section~3
of \cite{KP99}, and will become clear when we use them.

For example, if $C$ has precisely $r$ nodes, then its diagram consists
simply of $r$ roots of multiplicity 2; this diagram is denoted $r\IA_1$.
If $C$ has a simple cusp, then its diagram consists of three vertices: a
root of multiplicity 2, followed by a free vertex of multiplicity 1,
followed by a final vertex of multiplicity 1 and proximate to the root.
This diagram is denoted $\IA_2$.  Many more examples are discussed in
Section~2 of \cite{KP99}; in fact, there is there a classification of
all the $\ID$ with a single root $R$ and with $\cod(\ID)\le10$ (whence
$m_R\le 4$) and also of all those $\ID$ with $m_R\le3$.

The following lemma will be used to prove the first assertion of Theorem
(2.5).

 \lem4 Only finitely many distinct minimal Enriques diagrams arise from
the fibers of $D/Y$.
 \pf As $C$ ranges over the fibers, the numbers $\dim H^1(\cO_C)$ are
bounded, say by $p$.  Fix a relatively ample sheaf on $F/Y$; then the
numbers $\deg C$ are defined, and they too are bounded, say by $m$.  Fix
an arbitrary reduced $C$, and let $f\:C'\to C$ be the normalization map.
Then the number of connected components of $C$ is equal to $\dim
H^0(\cO_C)$, and the number of irreducible components of $C$ is equal to
$\dim H^0(\cO_{C'})$.  Hence
        $$\dim H^0(\cO_C)\ge 1 \and \dim H^0(\cO_{C'})\le m.$$

Consider the standard short exact sequence,
        $$0\to \cO_C\to f_*\cO_{C'}\to f_*\cO_{C'}/\cO_C\to0.$$
 In view of the preceding paragraph, this sequence yields the bound,
        $$\dim f_*\cO_{C'}/\cO_C \le p+m-1.$$

 Let $\ID$ be the diagram of $C$.  Then $\dim f_*\cO_{C'}/\cO_C =
\delta(\ID)$ by the Noether--Enriques theorem; see \cite{KP99,
Prop.~(3.1), p.~220}.  Hence, by the definition of $\delta(\ID)$,
          $$\sum_{V\in\ID,\ m_V\ge2}m_V/2\le\delta(\ID)\le p+m-1.$$
 Thus the number of vertices $V$ with $m_V\ge2$ is bounded, and the
weights $m_V$ themselves are bounded.

It remains to bound the number of vertices $V$ with $m_V=1$.  Each free
$V$ determines a distinct branch by \cite{KP99, Lem. (2.1), p.~214}.  So
the number of free $V$ is bounded by $r(\ID)$, which is equal to the
total number branches of $C$ through all of its singular points
\cite{KP99, Prop.~(3.1), p.~220}.  Hence, by Part~(i) of Lemma~(3.5) in
the next section, the number of $V$ with $m_V=1$ is bounded by $\sum
m_R$ where $R$ ranges over the roots of $\ID$.  However, $\sum m_R$ is
bounded by virtue of the last display in the preceding paragraph.  The
proof is now complete.\qed

We can now state our general theorem, and prove its first assertion.
The rest of the proof is found in \cite{K--P}.  (See also Section 4 of
\cite{KP99}.)

 \thm5 In the above setup, assume
 \rbit i  if\/ $Y(\infty)\ne\emptyset$, we have $\cod Y(\infty)\ge r+1$,
and
 \rbit ii for each $\ID$ such that $Y(\ID)\ne\emptyset$, we have
 $\cod Y(\ID)\ge\min(r+1,\,\cod\ID)$.\par
 \noindent Then either $Y(r\IA_1)$ is empty, or it has pure codimension
$r$; in either case, its closure $\?{Y(r\IA_1)}$ is the support of a
natural nonnegative cycle $U(r)$.  Furthermore, if $r\le8$, then the
rational equivalence class $[U(r)]$ is given by the formula
        $$[U(r)]=P_r(a_1,\ldots,a_r)/r! \where a_q:=\pi_*b_q$$
 and $b_q$ is a certain polynomial in $v,\,w_1,\,w_2$, namely, that
output by\/ \rm Algorithm (2.3).

\pf In the relative Hilbert scheme $\Hilb^r_{F/Y}$, form the open
subscheme $H(r)$ parameterizing the sets $G$ of $r$ distinct points in
the fibers of $F/Y$.  Re-embed $H(r)$ in $\Hilb^{3r}_{F/Y}$ by sending a
$G$ to the subscheme defined by the square of its ideal; this embedding
is well defined because $F/Y$ is smooth.  Next, form the intersection
	$$Z(r):=H(r)\bigcap \Hilb^{3r}_{D/Y},$$
 its closure $\?{Z(r)}$, and the fundamental cycle
$\bigl[\?{Z(r)}\bigr]$.  Push $\bigl[\?{Z(r)}\bigr]$ down to $Y$; the
result is, by definition, $U(r)$.

On $Y$, the image of $\?{Z(r)}$ contains, as a dense subset, the image of
$Z(r)$.  The latter image consists of $Y(\infty)$ plus the set of all
$y\in (Y-Y(\infty))$ such that $D_y$ has $r$ or more distinct singular
points.  The latter condition implies that the minimal Enriques diagram
$\ID$ of $D_y$ has $r$ roots or more.  So, $\cod(\ID)\ge r$, and
$\cod(\ID)=r$ if and only if $\ID=r\IA_1$ (because no final vertex, or
leaf, of $\ID$ can be a free vertex of multiplicity 1).  Hence, either
$y\in Y(\ID)$ with $\cod(\ID)> r$ or else $y\in Y(r\IA_1)$.  Also, the fiber
of $Z(r)$ over $y$ is finite, and it has cardinality 1 if $y\in Y(r\IA_1)$;
moreover, the image of $Z(r)$ contains $Y(r\IA_1)$.

Each component of $Z(r)$ is of dimension at least $\dim(Y)-r$, because
$H(r)$ is of dimension $\dim(Y)+2r$ and because $\Hilb^{3r}_{D/Y}$ is the
zero scheme of a regular section of a bundle of rank $3r$ on
$\Hilb^{3r}_{F/Y}$.  Now, by the preceding paragraph, the fibers of
$Z(r)/Y$ are finite off $Y(\infty)$, and have cardinality precisely 1 over
$Y(r\IA_1)$.  Moreover, the image of $Z(r)$ is contained in $Y(r\IA_1)$ plus
the union of $Y(\infty)$ and certain $Y(\ID)$ with $\cod(\ID)> r$; these
$\ID$ are finite in number by Lemma~(2.4) above.  Hence the hypotheses
of the theorem imply that $U(r)$ is a cycle of pure codimension $r$, and
its support is  $\?{Y(r\IA_1)}$.  The first assertion is now proved.\qed

Since the characteristic is 0, a map between integral schemes has
degree 1 if its fibers have cardinality one.  Hence the above
considerations also yield the following lemma, which we use in
conjunction with Theorem (2.5).

 \lem6 The enumerating cycle $U(r)$ is reduced if and only if the scheme
$Z(r)$ is reduced on an open set that dominates $Y(r\IA_1)$.
 \endproclaim\vskip-\bigskipamount

\rmk7 It is natural to conjecture that the theorem generalizes to any $r$.
More precisely, for any $r$, the hypotheses of the theorem should imply
that the class $[U(r)]$ is given by a universal polynomial in the classes
$y(a,b,c)$.  Moreover, this polynomial should be of the form
$P_r(a_1,\ldots,a_r)/r!$ where $a_q:=\pi_*b_q$ and $b_q$ is the
weighted homogeneous polynomial output by a suitable extension of
Algorithm (2.3) and evaluated at $v,\,w_1,\,w_2$.

It is also natural to conjecture that the theorem generalizes so as to
enumerate the $y\in Y$ such that $D_y$ has a given equisingularity type,
say that represented by a minimal Enriques diagram $\ID$.  More
precisely, set $r:=\cod\ID$, and let $\rho$ be the number of roots of
$\ID$.  Then the hypotheses should imply that the closure of $Y(\ID)$ is
the support of a natural positive cycle $U(\ID)$, and its class
$[U(\ID)]$ is given by a universal polynomial of degree $\rho$ in the
$y(a,b,c)$ with $a+b+2c\le r+2$.

Evidence for this conjecture is provided by the case of a suitably
general linear system on a fixed smooth irreducible projective surface.
First, Theorem (1.2) of \cite{KP99, p.~210}, enumerates the curves with
a triple point of a given type and additionally up to three nodes.
Second, this conjecture implies G\"ottsche's conjecture \cite{Gtt98,
Rmk.~5.4, p.~532}, which enumerates the curves with several ordinary
multiple points.

It is easy to construct a natural candidate for the cycle $U(\ID)$ by
generalizing the construction of $U(r)$.  Namely, set $d:=\deg(\ID)$,
and in $\Hilb^d_{F/Y}$ form the set $H_{F/Y}(\ID)$ of points
parameterizing the complete ideals with $\ID$ as associated diagram.
For example, $H(r\IA_1)=H(r)$.  (The set $H_{F/Y}(\ID)$ is studied in
Section~4 of \cite{KP99} and is studied further in \cite{K--P}.)  Owing
to the work of Nobile and Villamayor \cite{NV97, Thm.~2.6}, or that of
Lossen [pvt.\ comm.], $H_{F/Y}(\ID)$ is locally closed; in fact, it is a
smooth $Y$-scheme.  Form the intersection $Z(\ID):=H_{F/Y}(\ID)\bigcap
\Hilb^d_{D/Y}$, its closure $\?{Z(\ID)}$, and the fundamental cycle
$[\?{Z(\ID)}]$.  Push $[\?{Z(\ID)}]$ down to $Y$; take the result to be
$U(\ID)$.

The support of $U(\ID)$ contains, as a dense subset, the set of $y\in Y$
such that the Enriques diagram of $D_y$ contains $\ID$.  Hence the
hypotheses of the theorem imply that the support of $U(\ID)$ is equal to
the closure of $Y(\ID)$.

\sct3 Plane curves

Let $N_r(m)$ be the (unweighted) number of reduced plane curves of
degree $m$, that possess exactly $r$ (ordinary) nodes and that contain
$m(m+3)/2-r$ points in general position.  In this section, for $r\le8$,
we prove G\"ottsche's conjecture \cite{Gtt98, Conj.~4.1, p.~530}, about
$N_r(m)$; more precisely, we prove Theorem~(3.1), which is our first
main result.  The proof relies on Theorem~(2.5), which is solely
responsible for the restriction $r\le8$: if Theorem~(2.5) is proved for
more values of $r$, then Theorem~(3.1) will follow for these same
values.  We end the section with a survey of related work and with some
instructive examples.

\thm1 Assume $r\le8$ and  $m\ge r/2+1$.  Then
        $$N_r(m)=P_r(a_1,\ldots,a_r)/r!$$
 where $P_r$ is the Bell polynomial, defined by Identity \(2.2), and the
$a_q$ are the quadratic polynomials in $m$ listed in Table \Cs2).
\endproclaim
\goodbreak
\midinsert
$$\matrix{\hfil\hbox{\bf Table (3.2)}\hfil\cr
\hfil\hbox{\bf The polynomials $a_q(m)$ for plane curves }\hfil\strut\cr
\noalign{\medskip\hrule\medskip}
    a_1 = 3m^2-6m+3 = 3(m-1)^2\cr
    a_2 = -42m^2+117m-75 = -3(m-1)(14m-25)\cr
    a_3 = 1380m^2-4728m+3798\cr
    a_4 = -72360m^2+287010m-271242\cr
    a_5 = 5225472m^2-23175504m+24763752\cr
    a_6 = -481239360m^2+2334195360m-2748951000\cr
    a_7 = 53917151040m^2-281685755520m+359332109280\cr
    a_8 = -7118400139200m^2+39618359640720m-54066876993360\cr
 \noalign{\medskip\hrule\smallskip}}$$
 \endinsert

{\bf Proof.}  We apply Theorem~(2.5).  Let $Y$ be the
projective space parameterizing the plane curves of degree $m$, so $\dim
Y =m(m+3)/2$.  Set $S:=\IP^2$ and $F:=S\x Y$, and let $D\subset F$ be
the total space of curves.

Consider the set $Y(\infty)$ of $y\in Y$ such that the curve $D_y$ has an
$s$-fold component for some $s\ge2$.  If $m=1$, then $Y(\infty)$ is
empty.  Suppose $m\ge2$.  Then the $D_y$ with $s=2$ form a subset of
maximal dimension, namely, $(m-2)(m+1)/2+2$.  Hence $\cod Y(\infty)=
2m-1$.  Since $m\ge r/2+1$ by hypothesis, Hypothesis~(i) of Theorem~(2.5)
follows.  Furthermore, its Hypothesis (ii) holds owing to Parts (i) and
(ii) of Lemma~\Cs3) below.  Hence we may apply Theorem~(2.5).

Theorem~(2.5) implies that the closure of $Y(r\IA_1)$ is the support of
a nonnegative cycle $U(r)$, whose class is equal to
$P_r(a_1,\ldots,a_r)/r!\cdot h^r$ where the $a_q$ are certain integers
and where $h:=c_1(\cO_Y(1))$.  In fact, the argument at the top of
\p.232 of \cite{KP99} shows that the $a_q$ are equal to certain linear
combinations of the four basic Chern numbers $d$, $k$, $s$ and $x$.
These combinations are listed on \p.210 of \cite{KP99}.  Moreover, since
$S:=\IP^2$, the four numbers are, respectively, $m^2$, $-3m$, 9 and 3.
Formal calculations now yield the values in Table \Cs2).

Finally, $U(r)$ is reduced by Lemma~\Cs4) below.  Let $M\subset Y$ be the
linear space representing the plane curves that contain $m(m+3)/2-r$
points in general position.  Then $M\cap U(r)$ is finite, reduced, and
contained in $Y(r\IA_1)$ by Lemma (4.7) on \p.232 of \cite{KP99}.  Hence
$N_r(m)$ is equal to $P_r(a_1,\ldots,a_r)/r!$, and the proof is
complete.\qed

 \lem3 Assume $m\ge r/2 + 1$.  Let $Y$ be the projective space of plane
curves of degree $m$, and let $\ID$ be a (nonempty) minimal Enriques
diagram such that $Y(\ID)\ne\emptyset$.
 \rbit i If $\cod (\ID) \le r$, then $\cod(Y(\ID),Y) = \cod (\ID)$
and $Y(\ID)$ is smooth.  Moreover, then $Y(\ID)$ represents the functor
of $\ID$-equisingular families of plane curves of degree $m$ (their
parameter spaces need not be reduced).
 \rbit ii If $\cod(\ID)\ge r+1$, then $\cod(Y(\ID),Y) \ge r+1$.
 \pf
 Let $C$ be a curve corresponding to an arbitrary (closed) point of
$Y(\ID)$.  For a moment, suppose that $\ID$ consists of one vertex of
multiplicity $m$.  Then $\cod(\ID)={m+1\choose 2}-2$.  Furthermore, $C$
has an ordinary $m$-fold point.
 % (so $C$ is the union of $m$ concurrent lines).
  Hence $Y(\ID)$ is smooth, it represents
the functor, and $\cod(Y(\ID),Y) = \cod(\ID)$ owing to Greuel and
Lossen's \cite{GL96, Cor.~5.1 a), p.~339}.  Thus Parts (i) and (ii) hold
in this case.

For the rest of the proof, suppose therefore that $\ID$ does not
consists of one vertex of multiplicity $m$.  Then $C$ is not a union
of $m$ concurrent lines.  Now, $\deg C=m$.  So  $m\ge3$.

Let $\tau\es$ be the colength of the global equisingular ideal of
$C$ in $\IP^2$.  If $4m>4+\tau\es$, then $Y(\ID)$ is smooth, it
represents the functor, and $\cod(Y(\ID),Y) = \tau\es$ owing to Greuel
and Lossen's \cite{GL96, Cor.~3.9 b) and d), p.~334}, which applies
since $C$ is not a union of $m$ concurrent lines and since $m\ge3$.

Consider the multigerm of $C$ along its singular locus, a corresponding
miniversal deformation base space $B$, and the subspace of equisingular
deformations $B\es$.  Then $B\es$ is smooth and $\cod(B\es,B) = \tau\es$
owing to Wahl's \cite{Wahl74, Thm.~7.4, p.~162}.  However,
$\cod(B\es,B)=\cod(\ID)$ by \cite{KP99, Cor.~(3.3), p.~222}, (closely
related formulas were given by Wall \cite{Wall84, Thm.~8.1, p.~505}, by
Mattei \cite{Mat91, Thm.~(4.2.1), p.~323} and by T. de Jong
\cite{TdJ00, Thm.~3.5}; the present authors
are grateful to T. de Jong for pointing out the first two references).
Thus $\tau\es = \cod (\ID)$.

Suppose $\cod (\ID) \le r$.  Now, $ r\le 2m-2$ by hypothesis.  Also,
$m\ge2$; in fact, $m\ge3$.  So $2m-2<4m-4$.  Hence $4m-4>\cod (\ID)$.
So $4m>4+\tau\es$ by the preceding paragraph.  By the paragraph before
it, Part~(i) therefore holds.

Suppose that $\cod(\ID)\ge r+1$ instead.  If $4m-4>\cod (\ID)$, then as
in the preceding case, $\cod(Y(\ID),Y) = \cod (\ID)$, and so Part~(ii)
holds.  So suppose that $4m-4\le\cod (\ID)$.  Now, $\cod (\ID)\le
2\delta(\ID)$ by Part~(v) of Lemma (3.5).  Hence $2m-2\le\delta(\ID)$.
Now, $\delta(\ID)$ is equal to the genus discrepancy by the
Noether--Enriques theorem; see \cite{KP99, Prop.~(3.1), p.~220}.  Hence
$\cod Y(\ID) \ge \delta(\ID)$, and if equality holds, then $\ID
=\delta(\ID) \IA_1$, owing to Zariski's \cite{Z82, Thm.~2, p.~220}.
Now, for any $s$, we have $\cod(s \IA_1)=s$ and $\delta(s \IA_1)=s$.
Therefore, if $\cod Y(\ID) = \delta(\ID)$, then $\cod Y(\ID) =
\cod(\ID)$, and so Part~(ii) holds in this case.  However, if $\cod
Y(\ID) > \delta(\ID)$, then $ \cod Y(\ID)> 2m-2$ since
$2m-2\le\delta(\ID)$, and so Part~(ii) holds in any case.  The proof is
now complete.\qed

 \lem4 Consider the cycle  $U(r)$ of Theorem~(2.5).
If $m\ge r/2 + 1$,  then  $U(r)$  is reduced.

 \pf By definition, $U(r)$ is the image on $Y$ of the fundamental
cycle of the closure $\?{Z(r)}$ of the intersection
                   $Z(r):=H(r)\bigcap\Hilb^{3r}_{D/Y}$.
 By Lemma (2.6), $U(r)$ is reduced if $Z(r)$ is reduced on an open
set $Z^0$ that dominates $Y(r\IA_1)$.  We now construct such a
$Z^0$ by taking the inverse image of a suitable dense open subset
$Y(r\IA_1)^0$ of $Y(r\IA_1)$, and then we prove that the map
$Z^0\to Y$ factors through the reduced scheme $Y(r\IA_1)^0$ and that
the induced map $Z^0\to Y(r\IA_1)^0$ is an isomorphism.

Take any dense open subscheme $Y^0$ of $Y$ such that $Y^0\cap
\?{Y(r\IA_1)} \subset Y(r\IA_1)$, and denote the preimage of $Y^0$ in
$Z(r)$ by $Z^0$.  Taking $Y^0$ smaller if necessary, we may assume that
the map $Z^0\to Y^0$ is finite.  Set $D^0:=D\x_Y Z^0$.  Via the
projection to $H(r)$, view $Z^0$ as the parameter space of a flat family
of $r$ distinct points in the fibers of $F/Y$; denote the total space by
$W^0$.  Then, over a point of $Z^0$, the fiber of $W^0$ is just the set
of $r$ nodes of the the fiber of $D^0$.  Let $W_{\(2)}^0$ be the
infinitesimal thickening of $W^0$ defined by the square of its ideal.
Since $Z^0\subset Z(r)$, we have $W_{\(2)}^0\subset D^0$.

Let $\beta\:F^\star\to F\x_Y Z^0ñ$ be the blowup along $W^0$.  Set
$E^\star :=\beta^{-1}W^0$, so $E^\star$ is the exceptional divisor.  Set
$D^\star := \beta^{-1} D^0-2E^\star$.  Then $D^\star$ is effective since
$W_{\(2)}^0 \subset D^0$.  Moreover, the fibers of $D^\star/Z^0$ are the
proper transforms of the fibers of $D^0/Z^0$.  Hence $D^\star/Z^0$ is
smooth, and $(D^\star\cap E^\star)/Z^0$ is a family of $r$ pairs of
distinct points.  Thus, after localizing via the \'etale covering
$W^0/Z^0$, we obtain an $r\IA_1$-equisingular section of $D^0\x_{Z^0}
W^0/W^0$; in other words, $D^0/Z^0$ is an equisingular family of
$r$-nodal curves.  Now, $Y(r\IA_1)$ represents the functor of such
families by Part (i) of Lemma~(3.3).  Hence, the map $Z^0\to Y$ factors
through the reduced subscheme $Y(r\IA_1)$, so through its dense open
subscheme $Y^0\cap Y(r\IA_1)$.  Set $Y(r\IA_1)^0:=Y^0\cap Y(r\IA_1)$.

The map $Z^0\to Y(r\IA_1)^0$ is finite and surjective; moreover, its
fibers have cardinality 1 by the analysis in the middle of Section~2.
To prove that this map is an isomorphism, it suffices, since the
characteristic is 0, to prove that each closed fiber is reduced.
Suppose one isn't.  Then it contains a copy of Spec(A) where
$A:=\IC+\IC\epsilon$ is the ring of dual numbers.  Let $C$ be the
$r$-nodal curve in question.  Then $W^0\ox A/A$ is an \'etale family
supported on the set of nodes of $C$.  Furthermore, its infinitesimal
thickening $W_{\(2)}^0\ox A$ is contained in $C\ox A$.

This situation is untenable.  Indeed, work locally analytically at one
of the nodes of $C$.  Choose coordinates $X$, $Y$ so that $C:XY=0$.  Say
$W^0\ox A$ is defined by $X-a\epsilon=0$ and $Y-b\epsilon=0$.  Then the
ideal of $W_{\(2)}^0\ox A$ is generated by the three polynomials,
        $$X^2-2a\epsilon X,\ XY-\epsilon(aY+bX),\ Y^2-2b\epsilon Y.$$
 However, this ideal does not contain $XY$.  Thus the lemma is proved.\qed

 \lem5 Let $\ID$ be a minimal Enriques diagram with one root $R$.  Let
$\IS$ be the set of satellites of $\ID$, and set $e(\ID):=\mu(\ID)+m_R-1$.
Then
 \rbit i  $m_R=r(\ID)+\sum_{V\in \IS}m_V$;
 \rbit ii $\delta(\ID)\le \cod (\ID)$, with equality if and only if
$\ID=\IA_1$;
 \rbit iii $\cod (\ID) \le \mu(\ID)$, with equality if and only if
$\ID$ is $\IA_k$, $\ID_k$, $\IE_6$, $\IE_7$, or $\IE_8$;
 \rbit iv $\mu(\ID)\le 2\delta(\ID)$, with equality if and only if
$r(\ID)=1$;
 \rbit v $\cod (\ID) \le 2\delta (\ID)$, with equality if and only if
$\ID$ is either $\IA_{2i}$, $\IE_6$, or $\IE_8$;
 \rbit vi $2\delta(\ID)\le e(\ID)$, with equality if and only if
$m_R=r(\ID)$;
 \rbit {vii} $e(\ID)\le \cod (\ID)+\delta(\ID)$, with equality if
and only if $\ID=\IA_1$ or $\ID=\IA_2$;
 \rbit {viii} $e(\ID)\le 2\cod (\ID)$, with equality if and only if
$\ID=\IA_1$.
 \pf Consider Part~(i).  Let $V$ and $W$ be vertices, and write $W\succ_{\rm
imm} V$ if $W$ is an immediate successor of $V$.  If not, but $W$ is
proximate to $V$, then $W$ is a satellite of $V$.  By the ``law of
proximity,'' $W$ cannot also be a satellite of a second vertex.  Now, by
definition, $r(\ID):=\sum_V\bigl(m_V-\sum_{W\succ V}m_W\bigr)$.  Hence
        $$r(\ID)=\tsum_V\bigl(m_V-\sum_{W\succ_{\rm imm}
                V}m_W\bigr)-\sum_{V\in \IS}m_V.$$
 In the first sum, all the terms cancel except $m_R$.  Thus Part~(i)
holds.

Consider Part~(ii).  Denote the set of free vertices other than $R$ by
$\IF$.  Since $R$ is the only root, the definitions yield
  $$\eqalign{
\cod (\ID)-\delta(\ID)&=\textstyle{\sum_V}\bigl({{m_V+1}\choose 2}
                        -{m_V\choose2}\bigr)-1-\frs(\ID)\cr
        &=\textstyle{\sum_V} m_V-1-\frs(\ID)\cr
        &=(m_R-2)+\textstyle{\sum_{V\in \IF}} (m_V-1)
        +\sum_{V \in \IS} m_V.\cr}\eqno(3.5.1)$$
 Since $m_R\ge2$, the last term is nonnegative.  Moreover, it vanishes
if and only if $m_R=2$, and $m_V=1$ for all $V\in \IF$, and there are no
satellites.  However, the latter three conditions hold if and only if
$\ID=A_1$; see \cite{KP99, Section~2}.  Thus Part~(ii) holds.

Consider Part~(iii).  The definitions yield
   $$\mu(\ID)-\cod (\ID)=2\delta(\ID)-\deg(\ID)+2+\frs(\ID)-r(\ID).$$
Now, $\frs(\ID)$ is simply the total number of vertices less the number
of satellites; so
        $$\frs(\ID)= \tsum_V1 -\sum_{V\in \IS} 1. $$
 Hence Part (i) yields
  $$\eqalign{\mu(\ID)-\cod (\ID)
  &=\tsum_V\bigl(2{m_V\choose2}-{{m_V+1}\choose 2}+1\bigr)
        +(2-m_R)+\sum_{V\in \IS}(m_V-1)\cr
  &=\textstyle{m_R-2\choose2}
        +\sum_{V\in \IF}{m_V-1\choose 2}
        +\sum_{V\in \IS} {m_V\choose 2}.\cr}$$
 The last term is nonnegative. Moreover, it vanishes if and only if (1)
$m_R$ is 2 or 3, and (2) $m_V$ is 1 or 2 for all $V\in \IF$, and (3) $m_V$
is 1 for all satellites $V$.  However, the latter three conditions hold
if and only if $\ID$ is either $\IA_k$, $\ID_k$, $\IE_6$, $\IE_7$, or
$\IE_8$; see \cite{KP99, Section~2}.  Thus Part~(iii) holds.

Part~(iv) follows immediately from the definition of $\mu(\ID)$ since
$r(\ID)\ge1$.

Part~(v) follows immediately from Parts~(iii) and (iv) and from
Table~2-1 of \cite{KP99, p.~219}, which lists the values of $r(\ID)$ for
all the $\ID$ in question.

Consider Part~(vi).  The definitions of $e(\ID)$ and $\mu(\ID)$ yield
        $$e(\ID)-2\delta(\ID)=m_R-r(\ID).\eqno(3.5.2)$$
 Now, $m_R\ge r(\ID)$ by Part~(i), and Part~(vi) follows.

Consider Part~(vii).  Together (3.5.1) and (3.5.2) yield
 $$\cod (\ID)+\delta(\ID)-e(\ID)
        =r(\ID)-2+\tsum_{V\in \IF} (m_V-1)
        +\sum_{V\in \IS} m_V.$$
 Suppose $r(\ID)\ge 2$.  Then the right side is nonnegative.  Since
every free vertex of multiplicity 1 must be followed by a satellite, the
right side vanishes if and only if $r(\ID)=2$ and there are no other
vertices than the root, hence, if and only if $\ID=\IA_1$.

Suppose $r(\ID)=1$.  Then there is at least one satellite by Part~(i).
Hence, the right side is nonnegative.  It vanishes precisely when there
is just one free vertex other than $R$ and just one satellite, and both
have weight 1.  The latter condition implies $m_R=2$ by Part~(i) since
$r(\ID)=1$.  So the condition holds if and only if $\ID=\IA_2$.  Thus
Part~(vii) holds.

Finally, Part~(viii) follows immediately from Parts~(vii) and (ii).
Thus the lemma is proved.\qed

 \rmk6 Lemma~(3.5) is purely combinatorial.  However, it can be
interpreted geometrically, as asserting properties of an arbitrary curve
$C$ belonging to $Y(\ID)$.  Indeed, $r(\ID)$, $\delta(\ID)$, and the
other numerical characters of $\ID$ are equal to corresponding
characters of $C$.  All but $e(\ID)$ were treated in \cite{KP99, Section~3};
however, $e(\ID)$ is equal to the multiplicity of the Jacobian ideal,
owing directly to Teissier's work \cite{Te73, 1.6, p.~300}.

Correspondingly, Lemma~(3.5) can be proved via alternative geometric
arguments.  For example, the inequality $\cod (\ID) \le \mu(\ID)$ of
Part~(iii) just says that the modality $\mod(C)$ is nonnegative.
Indeed, $\mod(C)=\mu(C)-\tau\es$ by Greuel and Lossen's \cite{GL96,
Lem.~1.3, p.~326}, and $\tau\es = \cod (\ID)$ by the proof of
Lemma~(3.3) above.  Alternatively, $\mu(C)\ge\tau\es$ because the
equisingular ideal contains the Jacobian ideal by Wahl's work
\cite{Wahl74}; see the proof of Prop.~6.1, top of \p.159.

Similarly, the inequality $e(\ID)\le 2\cod (\ID)$ of Part~(viii) holds
because the equisingular, equiclassical, and equigeneric ideals form an
ascending chain.  Indeed, the inequality was proved this way by Greuel,
Lossen, and Shustin in \cite{GLS97, Lem.~2.2, p.~601}.  (However,
there's a typo in the proof: the colengths of the equiclassical and
equigeneric ideals are transposed.)  These authors and others write
`$\kappa$' instead of `$e$'.

 \rmk7
 The formula $N_1(m)=3m^2-6m+3$ was given by Steiner \cite{St48, p.~499},
in 1848, but it was probably known earlier.  After all, $N_1(m)$ is
simply the number of singular members of a general pencil of plane
curves of degree $m$.  So $N_1(m)$ is just the degree of the
discriminant of a general ternary form of degree $m$, viewed as a
polynomial in its coefficients, because the discriminant is the
resultant (or ``eliminant'') of the three partials.

The formula $N_2(m)=3/2(m-1)(m-2)(3m^2-3m-11)$ was given by Cayley
\cite{C66, Art.~33, p.~306}, in 1863.  He considered a varying pencil,
and formed its ``double discriminant.''  One factor has $N_2(m)$ as its
degree.  Cayley determined the degrees of the other two factors and the
degree of the double discriminant, then he divided.

The same formula was found a little differently by Salmon \cite{S65,
Appendix IV, p.~506}, (there is a typo: a `1' instead of an `11').
Salmon acknowledged Cayley's work, saying: ``Mr.~Cayley had arrived at
these numbers by a different process in a Memoir communicated to the
Cambridge Philosophical Society, but not yet published.''  First, Salmon
computed the number of curves with either two nodes or one cusp,
$9/2(m-1)(m-2)(m^2-m-1)$; then he subtracted the number of curves with
one cusp, $12(m-1)(m-2)$.  (See also the bottom of \p.361 in Salmon's
\cite{S79}.)

The formula for $N_2(m)$ was recovered implicitly via a third method by
S. Roberts \cite{R67, p.~276} in 1867.  At the bottom of \p.275, he said
his work agrees with Salmon's.

The formula $N_3(m)=9/2m^6-27m^5+9/2m^4+423/2m^3-229m^2-829/2m+525$ was
given implicitly by Roberts \cite{R75, pp.~111--12}, in 1875.  His
primary interest lay in his  new method for obtaining the degree of the
polynomial condition that three ternary equations have three common
solutions.  As an application, he discussed the theory of the
reciprocal, or dual, surface of a surface of degree $m$ in $\IP^3$.  In
effect, he determined the numbers $\beta$ of curves with one tacnode and
$\gamma$ of curves with one node and one cusp; he explicitly gave the
formulas,
        $$\eqalign{
        \beta&=50m^2-192m+168;\cr
        \gamma&=12(m-3)(3m^3-6m^2-11m+18).\cr}$$
 He also explicitly gave the formula,
        $$\beta+\gamma+N_3(m)
        =1/2(9m^6-54m^5+81m^4+63m^3-190m^2+11m+90).$$
 However, he did not solve for $N_3(m)$, which he denoted by $t$.

The formulas for $N_1(m)$, $N_2(m)$, and $N_3(m)$ were recovered
implicitly, and analogous formulas for $N_4(m)$, $N_5(m)$, and $N_6(m)$
were obtained explicitly, by Vainsencher \cite{V95, p.~515}, in 1995.
Vainsencher did not discuss the validity of these particular formulas,
but his general results, Propositions 3.5 and 4.1, do imply that there
exists some undetermined $m_0$ such that, for $m> m_0$, the formulas do
hold.

The formulas for $N_1(m)$, $N_2(m)$, and $N_3(m)$ were recovered
explicitly by Harris and Pandharipande \cite{HP95} later in 1995.  They
used an interesting new method, involving the geometry of the Hilbert
scheme of points in $\IP^2$ and the Bott residue formula.  In 1997, in
the paragraph before Definition 3.4 and in Prop.~3.6 in \cite{Ch97},
Choi derived from Ran's Theorem~5 in \cite{Ran89} that $N_r(m)$ is, for
every $r$ and $m>r$, given by a polynomial in $m$ of degree $2r$.

 A recursive formula was obtained by Caporaso and Harris \cite{CH98,
p.~353}, in 1998, which makes it possible to compute $N_r(m)$ for every
$r$ and $m$.  From this formula, though, it is not at all clear that,
when $r$ is fixed, $N_r(m)$ is given by a polynomial of degree $2r$ in
$m$ for $m> m_0$ for some $m_0$.  Nevertheless, as G\"ottsche observed
in \cite{Gtt98, Rmk.~4.2, p.~530}, and Choi observed in \cite{Ch97,
Rmk.~3.5.2}, if it is assumed that $N_r(m)$ is given by such a ``node''
polynomial for a known $m_0$, then it is possible to work out the
coefficients.

Given a specific value for $m_0$, such as Choi's value $m_0=r$ mentioned
above or the value $m_0=3r-1$ for $r\le8$ given in Thm.~(1.3) of
\cite{KP99}, it is possible to use Caporaso and Harris's formula to
check the validity of the values given by the polynomial for $m_0\ge
m\ge r/2+1$.  Thus, given an $r$ and an $m_0$, it is possible to prove
G\"ottsche's conjecture \cite{Gtt98, Conj.~4.1}, and so, for $r\le8$, to
obtain another proof of Theorem~(3.1).

\eg8 It is useful to look at some basic examples.  First, note that, for
any given $m$, there are several special ranges for $r$.  For
$r\le\min(2m-2,\,8)$, the formula $N_r(m)=P_r(m)/r!$ holds by
Theorem~(3.1).  For $r=2m-1$, the sets $Y(r\IA_1)$ and $Y(\infty)$ have
the same dimension when both are nonempty; see the first part of the
proof of that theorem.  Both sets are empty when $m=1$, but $Y(\infty)$
is nonempty for $m\ge2$.  For $r>(m-1)(m-2)/2$, if $y\in Y(r\IA_1)$,
then $D_y$ is reducible; otherwise, $D_y$ would have strictly negative
geometric genus since $(m-1)(m-2)/2$ is equal to its arithmetic genus
$p$.  For $r=m(m-1)/2$, if $y\in Y(r\IA_1)$, then $D_y$ is the union of
$m$ lines.  Finally, for $r>m(m-1)/2$, the set $Y(r\IA_1)$ is empty;
indeed, applying the argument in the middle of Lemma~(2.4) with
$C:=D_y$, we see that $r\le p+m-1$, with equality if and only if $D_y$
has $m$ irreducible components.

Suppose $m=1$.  Then $D_y$ is a line for every $y\in Y$.  So $N_0(1)=1$,
and $N_r(1)=0$ for $r\ge1$.  On the other hand, direct computation
yields $P_0(1)=1$, and $P_r(1)=0$ for $r=1,2$, but $P_3(1)/3!=75$.  Thus
$N_r(1)=P_r(1)/r!$ holds for $r=0,1,2$ but not for $r=3$.  The
hypothesis $m\ge r/2+1$ of Theorem~(3.1) fails for $r\ge1$.  However,
the hypotheses of Theorem~(2.5) are satisfied for every $r$.  Hence, the
proof of Theorem~(3.1) shows that the formula $N_r(1)=P_r(1)/r!$ must
hold for $r=0,1,2$.  For $r\ge3$, the value of $P_r(1)$ is irrelevant
since $U(r)$ vanishes by reason of dimension.

Suppose $m=2$.  Then $N_r(2)=P_r(2)/r!$ holds for $r=0,1,2$ by
Theorem~(3.1).  If $y\in Y(\IA_1)$, then $D_y$ is a line-pair.  So
$N_1(2)$ is the number of line-pairs through four points in general
position.  This number is ${4\choose2}/2$, or 3, since two of the four
points will determine one of the lines, and the remaining two point will
determine the other.  Now, $Y(\ID)$ is empty for any $\ID$ other than
$\IA_1$.  So $N_r(2)=0$ for $r\ge2$.  On the other hand, direct
computation yields $P_1(2)=3$ and $P_2(2)=0$, but $P_3(2)/3!=-32$.  Thus
$N_r(2)=P_r(2)/r!$ checks for $r=1,2$.

The equation $N_r(2)=P_r(2)/r!$ fails, however, for $r=3$, although
Theorem~(2.5) nearly applies.  Indeed, all the relevant $Y(\ID)$ are
empty, and $Y(\infty)$ has its expected codimension, namely 3, but
Hypothesis~(i) requires $\cod Y(\infty)>3$.  In fact, $Y(\infty)$ is the
Veronese surface in $Y=P^5$.  So two general hyperplanes intersect each
other and $Y(\infty)$ in four distinct points.  If each hyperplane
represents the conics that contain a given point in $\IP^2$, then the
four points of intersection coalesce in the point that represents the
double-line through the two points in $\IP^2$.  Since $P_3(2)/3!=-32$,
this double-line may be interpreted as four coincident double-lines,
each the equivalent of $-8$ three-nodal conics.

Finally, consider the case $m=5$ and $r=8$.  Here, $N_8(5)$ is the
number of 8-nodal quintics through 12 points in general position.  Since
$8>(5-1)(5-2)/2$, these quintics are reducible.  So each is either the
union of two smooth conics and a line, or the union of a nodal cubic and
a line-pair.  Hence
  $$N_8(5)=\textstyle{12\choose 5}{7\choose 5}\big/2
  +N_1(3){12\choose 8}{4\choose2}\big/2 = 8316+17820=26136. $$
 Therefore, Theorem~(3.1) implies that $P_8(5)=26136\times8!\/$.

The value of $P_8(5)$ can be used to determine the multiplier $3281\cdot
7!\,$ of $x_4$ in $b_8$ in Algorithm (2.3).  Indeed, as indicated in
\cite{KP99, Section~4}, residual intersection theory shows that $x_4$ appears
with some multiplier, whose value does not vary with $F/Y$ and $D$.
Hence this value can be determined from any particular example that can
be worked out by other means.

\sct4 Threefolds in four-space

In this section, we enumerate the 6-nodal plane curves on a general
threefold in 4-space, recovering Vainsencher's formula.  In fact, we
correct a typo: the multiplier 5 appearing in Theorem~\Cs1) is lacking
in \cite{V95, \p.522}.  More importantly, we establish the formula's
validity, which was left unaddressed in \cite{V95}.  Then, in
Theorem~\Cs3), we rederive the number of {\it irreducible} 6-nodal plane
quintic curves on a general quintic threefold; again we follow
Vainsencher's approach \cite{V95, \pp.523--24}, (recovering his number),
but also establish the number's validity.  Its significance is recalled
in the introduction.

 \thm1 In $\IP^4$, consider a general threefold $Q$ of degree $m$.  If
$m\ge4$, then $Q$ contains precisely the following number of $6$-nodal
plane curves of degree $m$:
 $$\eqalign{
 &5(m^{18}-12m^{17}+24m^{16}+
  155m^{15}-405m^{14}+1082m^{13}-18469m^{12}\hfill\cr
 &\qquad+66446m^{11}-192307m^{10}+
  1242535m^9-4049006m^8+11129818m^7\cr
 &\qquad-53664614m^6+166756120m^5-415820104m^4+1293514896m^3\cr
 &\qquad\qquad-2517392160m^2 +1781049600m)/6!\,.\cr}$$
 \pf
  Let's apply Theorem (2.5) again.  Fix $r\le6$.  Let $Y$ be the
Grassmann variety of $2$-planes $H$ in $\IP^4$, and $\cQ$ the
tautological rank-$3$ quotient of $\cO_Y^5$.  Let
                    $$F:=\IP(\cQ)\subset Y\x \IP^4\eqno(4.1.1)$$
 be the total space of planes, $\pi\:F\to Y$ and $p\: F\to \IP^4$ the
projections.  Set
                         $$D:=F\cap (Y\x Q).$$

Since $Q$ is smooth, $Q$ contains no 2-plane $H$; otherwise,
there'd be a normal-bundle surjection, $\cN_{H/\IP^4}\onto
\cN_{Q/\IP^4}|H$, in other words, $\cO_H(1)^2\onto \cO_H(m)$; so
$3\times 2\ge (m+2)(m+1)/2$, contradicting $m\ge3$.  Since $Q$ contains
no $H$, each $D_y$ is a curve of degree $m$.  Hence $D$ is a relative
effective divisor on $F/Y$.  Moreover,
              $$\cO_F(D)=p^*\cO_{\IP^4}(m).\eqno(4.1.2)$$

Consider the set $Y(\infty)$ of $y\in Y$ such that the curve $D_y$ has a
multiple component.  To check Hypothesis (i) of Theorem (2.5), it
suffices to prove  $\cod Y(\infty)= 2m-1$ if $Y(\infty)\ne\emptyset$,
because $2m-1>r$ as $m\ge4$ and $r\le6$.  Let's count constants.

Let $U$ be the scheme of all smooth quintic threefolds, and $P$ the
scheme of all quintic curves in the fibers of $F/Y$ (so $P$ is
$\IP(\cSym^5(\cQ^*)$).  Form the natural map,
         $$\lambda\:U\x Y\to P;$$
 it's given by $\lambda(Q,H) = Q\cap H$.  Now, $\lambda$ is the
restriction to $U\x Y$ of a family over $Y$ of linear projections; each
has, as domain, the projective space of all quintic threefolds, and as
center at $H\in Y$, the subspace of threefolds containing the plane $H$.
Hence, $\lambda$ is smooth, and has irreducible fibers (conceivably
$\lambda$ is not surjective).

Let $P_\infty\subset P$ be the subset of quintic curves with a multiple
component.  The fibers of $P_\infty/Y$ have codimension $2m-1$ by the
analysis at the beginning of the proof of Theorem~(3.1).  So,
$\lambda^{-1}P_\infty$ has codimension $2m-1$ if it's nonempty.  Now,
$Y(\infty)$ is just the fiber of $\lambda^{-1}P_\infty$ over $Q\in U$,
and $Q$ is general.  Hence $\cod Y(\infty)= 2m-1$ if $Y(\infty)\ne
\emptyset$, as desired.  Thus Hypothesis (i) of Theorem (2.5) is
satisfied.

Hypothesis (ii) is also satisfied, as the same argument shows.  Indeed,
in the plane, the corresponding condition is satisfied by Lemma~(3.3),
which applies since $m\ge6/2+1$.

Therefore, we may apply Theorem (2.5).  We conclude that the closure
$\?{Y(r\IA_1)}$ is the support of a nonnegative cycle $U(r)$, whose class
$[U(r)]$ is given by the formula
        $$[U(r)]=P_r(a_1,\ldots,a_r)/r! \where a_q:=\pi_*b_q$$
 and $b_q$ is the polynomial in $v,\,w_1,\,w_2$ output by\/ \rm
Algorithm (2.3).

In the present case, $U(r)$ is reduced, as the same argument shows when
developed a little further.  Indeed, in the plane, the corresponding
cycle is reduced by Lemma~(3.4).  Moreover, by Lemma~(2.6), $U(r)$ is
reduced if and only if $Z(r)$ is reduced on an open set that dominates
$Y(r\IA_1)$.  So consider the open subset on which $Z(r)$ is reduced.
This subset can be shown to dominate $Y(r\IA_1)$ by developing the same
argument as the formation of $Z(r)$ and $Y(r\IA_1)$ is compatible with
the constructions involved.  Therefore, the desired number of curves is
equal to the degree $\int[U(6)]$.

 To compute the degree $\int[U(6)]$, set
        $$h:=c_1(\cO_{\IP^4}(1)) \and  q_i:=c_i(\cQ).$$
 Recall the Euler exact sequence:
        $$0\to\Omega^1_{F/Y}(1)\to\pi^*\cQ \to\cO_F(1)\to0.$$
 It and Formula (4.1.2) yield the following formulas:
        $$v=mp^*h\and w_1=\pi^*q_1 - 3p^*h \and w_2=\pi^*q_2 -2p^*h\pi^*q_1
        +3p^*h^2.\eqno(4.1.3)$$
 By reason of dimension, $h^j=0$ for $j>4$, and $\pi_*p^*h^j=0$ for
$j=0,1$; moreover, it is well known and easy to see that
  $$\pi_*p^*h^2=[Y]\and\pi_*p^*h^3=q_1\and\pi_*p^*h^4=q_1^2-q_2.
        \eqno(4.1.4)$$
 It is now a mechanical matter to compute the $a_q$, and then $[U(6)]$, as
polynomials in $m$, $q_1$, and $q_2$.  Finally, standard Schubert
calculus on $Y$ yields the following degrees:
 $$\int q_1^6=5,\ \int q_1^4q_2=3,\ \int q_1^2q_2^2=2,\ \int q_2^3=1.
        \eqno(4.1.5)$$
 It is now a mechanical matter to complete the proof.\qed

\rmk2 Other enumerations fall out of the proof of Theorem (4.1).
Indeed, the proof shows that the cycle $U(r)$ is reduced for every $r$,
and that its class can be computed mechanically as a polynomial in $m$,
$q_1$, and $q_2$.  It follows, for example, that in $\IP^4$ a general
threefold of degree $m\ge4$ contains precisely
        $$5m^9/6-5m^8+11m^7/2+23m^6/6+17m^5+359m^4/6
        -1165m^3/3+1024m^2/3+40m$$
 3-nodal plane curves whose plane meets three lines in general position.

Indeed, these curves are enumerated by the intersection of $U(3)$ and the
three special Schubert cycles defined by the three lines.  The intersection
is reduced by the theorem of transversality of the general translate
\cite{Kl74, Cor.~4, p.~291}, (which applies because the characteristic is
0).  The class of each Schubert cycle is $q_1$.  Hence the number of
curves is just $\int[U(3)]q_1^3$, and its value can be computed
mechanically.

\thm3 In $\IP^4$, a general quintic threefold contains precisely
$17,601,000$ {\rm irreducible} $6$-nodal plane curves of degree $5$.
 \pf (Compare with Vainsencher \cite{V95, \pp.523--24}.)  For $m=5$, the
formula in Theorem (4.1) yields the number 21617125.  From it, we must
subtract the number of reducible curves.  Each is, plainly, one of two
types: (1) the transversal union of a smooth conic and a smooth cubic,
or (2) the transversal union of a line and a binodal quartic.  Each
reducible curve is counted with multiplicity 1; indeed, all 21617125
curves are by the theorem.  So we may count set-theoretically.

A general quintic $Q$ contains precisely 609250 smooth conics, thanks to
the work of S. Katz (see \cite{Ka92, (2), p.~175}).  Each conic
determines a plane, and it meets $Q$ partially in the conic and
residually in a cubic.  The cubic is smooth and meets the conic
transversally because $Q$ is generic.  Indeed, form the space of all
triples $(H,A,B)$ where $H$ is a plane, $A$ is a conic in $H$, and $B$
is a cubic in $H$.  Form the subset $U$ of $(H,A,B)$ such that $A$ and
$B$ are smooth and meet transversally in $H$.  Clearly $U$ is open and
dense.  So its complement $R$ has smaller dimension.  Hence, counting
constants as in the proof of Theorem (4.1), we conclude that there is no
$(H,A,B)\in R$ such that $A\cup B\subset Q$.  Thus the first subtrahend
is 609250.

A general quintic $Q$ contains precisely 2875 lines (see \cite{Ka92,
(1), p.~175}).  Fix one, $L$ say.  We must enumerate those planes $H$
that meet $Q$ partially in $L$ and residually in a binodal quartic that
meets $L$ transversally.  We do so by building on the proof of
Theorem~(4.1); in particular, we use its notation.

In the Grassmannian $Y$ of all planes $H$, form the Schubert variety
$Y_L$ of all $H$ that contain $L$.  Correspondingly, in the total space
$F$, form the preimage $F_L:=\pi^{-1}Y_L$.  Also, in the space $U$ of
all smooth quintic threefolds, form the subspace $U_L$ of those that
contain $L$.  Then $Q\in U_L$.  In fact, we may assume that $(Q,L)$ is a
general pair in the space of all such pairs, for this space is
irreducible by Katz's Lemma 1.4 of \cite{Ka86, p.~153}.  In particular,
$Q$ is a general point of $U_L$.  Finally, set
        $$D_L:=F_L\cap (Y_L\x Q)\and D':=D_L-(Y_L\x L);\eqno(4.3.1)$$
 so the fibers of $D_L/Y_L$ are the quintic curves cut out by $Q$, and
the fibers of $D'/Y_L$ are the residual quartic curves.

In the space $P$ of all quintic curves in the fibers of $F/Y$, form the
subspace $P_L$ of those with $L$ as a component, and form the natural
map,
         $$\lambda_L\:U_L\x Y_L\to P_L,$$
 the restriction of $\lambda$.  So $\lambda_L$ is the restriction to
$U_L\x Y_L$ of a family over $Y_L$ of linear projections; each has, as
domain, the projective space of all quintic threefolds containing $L$.
In the latter projective space, $U_L$ is open and nonempty.  Hence,
$\lambda_L$ is smooth, and has irreducible fibers.

Let's apply Theorem (2.5) with $r:=2$ to $D'/Y_L$.  Hypotheses (i) and
(ii) can be checked by counting constants just as in the proof of
Theorem~(4.1), after replacing $Y$, $U$ and $P$ with $Y_L$, $U_L$ and
$P_L$, because $P_L$ may be viewed as the space of all quartics in the
fibers of $F_L/Y_L$.  By the same token, the cycle of binodal quartics
$D'_y$ is reduced.  Similarly, each $D'_y$ meets $L$ transversally,
because, given a plane $H$ containing $L$, in the space of binodal
quartics in $H$, those tangent to $L$ form a closed subset of
codimension 1, as yet another count of constants shows.  Therefore, we
may apply Theorem (2.5), and use it as follows to get the multiplier
of 2875.

We must find the classes $v':=[D']$ and $w_i':=
c_i(\Omega^1_{F_L/Y_L})$.  First, $w_i'=w_i|F_L$ since
$F_L:=\pi^{-1}Y_L$.  Now, $F:=\IP(\cQ)$ by (4.1.1), so
$F_L=\IP(\cQ|Y_L)$.  Hence, the inclusion of $Y_L\x L$ in $F_L$
corresponds to a surjection $\cQ|Y_L\onto \cO_{Y_L}^2$.  Denote its
kernel by $\cK$.  Then there is a natural exact sequence,
 $$\cK\ox\cSym(\cQ|Y_L)[-1]\to\cSym(\cQ|Y_L)\to\cSym(\cO_{Y_L}^2)\to0.$$
Passing to associated sheaves, we see that the image of
$(\pi_L^*\cK)(-1)$ in $\cO_{F_L}$ is just the ideal of $L$, where
$\pi_L\:F_L\to Y_L$ is the structure map.  Hence
        $[L] = -c_1\bigl((\pi_L^*\cK)(-1)\bigr)$.
 Now, we have $c_1(\cK)=c_1(\cQ|Y_L)=q_1|Y_L$.  Hence
$[L]=(p^*h-\pi^*q_1)|F_L$.  Therefore, (4.3.1) and (4.1.3) yield
$v'=(4p^*h+\pi^*q_1)|F_L$.

Finally, by Theorem (2.5) and the projection formula, the desired
multiplier is equal to $(1/2)\deg P_2(a_1,a_2)$ where, by the projection
formula,
         $$a_q= [Y_L]\cdot\pi_*b_q(4p^*h+\pi^*q_1, w_1,w_2).$$
 Now, by standard Schubert calculus, $[Y_L]=(q_1^2-q_2)^2$.  Owing to
(4.1.3), to (4.1.4) and to (4.1.5), it is now a mechanical matter to
compute the multiplier; its value is 1185.  Therefore, the number of
6-nodal degree $5$ plane curves on $Q$ is
        $$21617125-609250-2875\x1185=17601000,$$
 and the proof is complete.\qed

\sct5 Abelian surfaces

Fix an Abelian surface $A$.  In this section, assuming certain
genericity conditions, we enumerate the reduced and irreducible
curves $C\subset A$ satisfying these three conditions: they lie in a
given algebraic homology class with positive self-intersection, $d$ say;
they have given geometric genus $g$; and they pass through the
appropriate number of general points.  The appropriate number is $g$, as
we see while proving Theorem~(5.2), our main result of this section.

Each such curve $C$ must have a certain number of singular points, $r$
say; we prove that, under our genericity conditions, all $r$ points
are ordinary nodes.  Since $g+r$ is the arithmetic genus of $C$, and
since the canonical class of $A$ is trivial, the numbers $d$, $g$, and
$r$ are, owing to the adjunction formula, related by the equation
        $$d=2g+2r-2.$$

 The number of these $C$ turns out to depend only on $r$ and $g$, and
not on $A$; here we must also assume $r\le 8$.  So, let's denote the
number of $C$, as Bryan and Leung \cite{BL98} do, by $N_{g,r}$.  (They
impose no bound on $r$, but do require $A$ to be generic among the
Abelian surfaces for which the given homology class is algebraic; whence
this class must generate $\Pic(A)/\Pic^0(A)$.)  For $r$ fixed, $N_{g,r}$
is given by a ``node'' polynomial in $g$ of degree $r+1$.  The
polynomials are presented in Table~(5.1).
 \midinsert
 \centerline{\bf Table (5.1)}
 \centerline{\bf Node polynomials for $N_{g,r}$ \strut}
\kern-0.75\baselineskip
% \kern-1.5\normalbaselineskip
$$  \eightpoint
  \setbox0=\hbox{${}=24g(g-1)(81g^6-1053g^5+7200g^4-29970g^3+75814g^2
        -106347g+62685)/35$}
\eqalign{\noalign{\hrule \medskip}
N_{g,0}&=g\cr
N_{g,1}&=6g(g-1)\cr
N_{g,2}&=6g(g-1)(3g-4)\cr
N_{g,3}&=4g(g-1)(9g^2-27g+25)\cr
N_{g,4}&=6g(g-1)(9g^3-45g^2+94g-75)\cr
N_{g,5}&=12g(g-1)(27g^4-198g^3+687g^2-1213g+860)/5\cr
N_{g,6}&=4g(g-1)(81g^5-810g^4+4095g^3-11835g^2+18409g-11800)/5\cr
N_{g,7}&=24g(g-1)(81g^6-1053g^5+7200g^4-29970g^3+75814g^2
        -106347g+62685)/35\cr
N_{g,8}&=3g(g-1)(486g^7-7938g^6+69930g^5-389970g^4+1413384g^3\cr
 &\kern\wd0
 \llap{$-3216332g^2+4143290g-2279375)/35$}\cr
 \noalign{\medskip\hrule}}$$\endinsert
 Our genericity conditions are specified in the following theorem.

\thm2 In the above setup, assume   $A$ has Picard number $1$, and
say the given homology class is $m$ times the positive primitive class.
Assume either
 \rbit i that $m=1$ and $g> 5r+7$ or\par
 \rbit ii that $m\ge2$ and $g>(3m^2r+3m^2-2mr+2m+2r-2)/(2m-2)$.\smallskip
 \noindent Then the polynomial formulas in Table~(5.1) are valid.

\pf Yet again, we apply Theorem~(2.5).  Let $Y$ be the connected
component of $C$ in the Hilbert scheme of $A$.  Then $Y$ parameterizes
the curves algebraically equivalent to $C$; hence, $Y$ parameterizes the
curves homologically equivalent to $C$, since an Abelian variety has no
torsion (see \cite{B--L99, Prop.~7.1, p.~59}).  Set $F:=A\x Y$, and let
$D\subset F$ be the universal divisor.  To handle $Y$ and $D$, we need a
well-known description of them, which we now recall.

Fix an invertible sheaf $\cL$ on $A$ representing the given
homology class.  Denote by $\wh A$ the dual abelian surface, and by
$\cP$ the Poincar\'e bundle on $A\x \wh A$, which is trivial along
$0\x\wh A$ and $A\x0$.
On $\wh A$, form the direct image
        $$\cQ:={p_2}_*(\cP\ox p_1^*\cL).\eqno(5.2.1)$$
 If $\cN$ is a fiber of $\cP\ox p_1^*\cL$, then $H^q(\cN)=0$ for $q\ge1$
by the Kodaira vanishing theorem, because $\cN$ is algebraically
equivalent to $\cL$, so ample, and because the canonical bundle of $A$
is trivial.  Hence $\cQ$ is locally free, and its formation commutes
with base change.

  Because of this commutativity, the rank of $\cQ$ may
be determined on the fibers, where we may use the Riemann--Roch theorem
and the vanishing of the higher cohomology groups; thus,
        $$\rk(\cQ) = \int c_1(\cL)^2/2=d/2.\eqno(5.2.2)$$
 Then $Y:=\IP (\cQ^*)$ where $\cQ^*$ is the dual.  Hence $Y$ is smooth
and irreducible, and
  $$\dim Y=\rk Q -1+2=d/2+1=g+r.$$
 Thus $g$ is the appropriate number, as asserted in the first paragraph
of the section.

To construct $D\subset F:=A\x Y$, form the natural Cartesian diagram:
 $$\CD
  F @>1\x p>> A\x\wh A \\
 @VV\pi V          @VV p_2 V \\
    Y      @>p>>    \wh A \rlap{\thinspace.}
  \endCD$$
 On $Y$, the tautological map $p^*\cQ^*\to \cO_Y(1)$ and the base-change
isomorphism induce the composition
        $$\cO_Y(-1)\to p^*{p_2}_*(\cP\ox p_1^*\cL)\risom
   \pi_*(1\x p)^*(\cP\ox p_1^*\cL).$$
 Form its adjoint on $F$,
  $$\alpha\:\pi^*\cO_Y(-1)\to (1\x p)^*(\cP\ox p_1^*\cL).$$
 The zero locus of $\alpha$ is the universal divisor $D/Y$.  Therefore,
  $$\cO_{F}(D)=(1\x p)^*\cP\ox p_1^*\cL\ox \pi^*\cO_Y(1)\eqno(5.2.3)$$
 where now $p_1\:F\to A$ is the projection.

In order to apply Theorem (2.5), we must check its hypotheses, (i) and
(ii).  Now, each of the hypothesis, (i) and (ii), of Theorem~(5.2)
implies, owing to Lemma~(5.3) below, that, if $\cN$ is any fiber of
$\cP\ox p_1^*\cL$, then $\cN$ is $k$-very ample for $k:=3(r+1)-1$.  In
other words, given any ideal $\cI\subset \cO_A$ of colength $3(r+1)$ or
less, the natural map
        $$H^0(\cN)\to H^0(\cN/\cI\cN)\eqno(5.2.4)$$
 is surjective.

Let $\ID$ be a minimal Enriques diagram, and set $e:=\deg\ID$ and
$s:=\cod\ID$.  We must bound $\cod Y(\ID)$ in terms of $s$.  Form the
subset $H(\ID)$ of $\Hilb^e_A$ of complete ideals with $\ID$ as diagram.
Then $H(\ID)$ is locally closed, smooth, and equidimensional of
dimension $e-s$ by Proposition~(3.6) of \cite{KP99, p.~225}.  Furthermore,
$Y\x H(\ID)$ is equal to the subset $H_{F/Y}(\ID)$ of $\Hilb^e_{F/Y}$,
which was discussed in Remark~(2.7).  As in that remark, set
$Z(\ID):=H_{F/Y}(\ID)\bigcap \Hilb^e_{D/Y}$.

Suppose $e\le 3(r+1)$.  Let $\cI$ be an arbitrary complete ideal with
diagram $\ID$.  Then $\cI$ has colength $e$ by \cite{KP99, (3.4),
p.~223}.  So, if $\cN$ is any fiber of $\cP\ox p_1^*\cL$, then (5.2.4)
is surjective.  Hence
  $$\dim H^0(\cI\cN)= \dim H^0(\cN)-e.$$
 Furthermore, $H^1(\cN)=0$, as noted above.  Hence $H^1(\cI\cN)=0$.

Let $\cI^\dagger$ be the universal ideal on $A\x H(\ID)$, and
form these two sheaves:
        $$\cF:=p_{13}^*\cI^\dagger\cdot p_{12}^*(\cP\ox p_1^*\cL)
 \and \cR:={p_{23}}_*\cF$$
 where the $p_{ij}$ are the projections from $A\x\wh A\x H(\ID)$.  Every
fiber of $\cF$ is of the form $\cI\cN$, and as just noted,
$H^1(\cI\cN)=0$.  Hence $\cR$ is locally free on $\wh A\x H(\ID)$, and
the formation of $\cR$ commutes with base change.  Furthermore, $\rk\cR
=\dim H^0(\cI\cN)$.

We have $\IP(\cR^*)=Z(\ID)$ in $Y\x H(\ID)$.  Indeed, the inclusion
$\cF\to p_{12}^*(\cP\ox p_1^*\cL)$ and the base-change isomorphism
induce a composition on $\wh A\x H(\ID)$,
        $$\cR\to {p_{23}}_*p_{12}^*(\cP\ox p_1^*\cL)\risom p_1^*\cQ,$$
 where the $p_1$'s are different first projections.  Dualizing gives a
map, $p_1^*\cQ^*\to\cR^*$.  It is surjective because its formation
commutes with base change and because its fibers are surjective, since
each is the dual of an inclusion of vector spaces of the form
$H^0(\cI\cN) \to H^0(\cN)$.  Thus $\IP(\cR^*)$ naturally embeds in
$\IP(p_1^*\cQ^*)$, or $Y\x H(\ID)$.

A nonzero map $\sigma\:\cO_A\to\cN$ factors through $\cI\cN$ if and only
if the corresponding map $\sigma\ox\cN^*\:\cN^*\to\cO_A$ factors through
$\cI$, so if and only if the divisor defined by $\sigma$ contains the
finite subscheme defined by $\cI$.  Thus $\IP(\cR^*)$ and $Z(\ID)$ have
the same sets of closed points.  An analogous argument shows that they
have the same sets of $T$-points for any $T$; whence, as schemes,
$\IP(\cR^*)=Z(\ID)$.

Therefore, $Z(\ID)$ is a smooth and equidimensional scheme since $H(\ID)$ is
so, and
        $$\eqalign{\dim Z(\ID)
        &= \rk\cR -1 + \dim (H(\ID)\x \wh A)\cr
        &= \rk\cQ -e -1 + e-s +\dim\wh A\cr
        &= \dim(Y)-s.\cr}$$
 So, if the image of $Z(\ID)$ in $Y$ contains a nonempty set $S$, then
$\cod S\ge s$.

To check the hypotheses of Theorem (2.5), we apply the preceding
conclusion about $S$ in three cases.  First, take $\ID:=(r+1)\IA_1$.
Then $e=3(r+1)$ and $s=r+1$.  Furthermore, the image of $Z(\ID)$ in $Y$
contains $Y(\infty)$.  Hence, $\cod Y(\infty)>r$ if $Y(\infty)\ne
\emptyset$.  Thus Hypothesis (i) of Theorem (2.5) holds.

Second, let $\ID'$ be a diagram with $Y(\ID')\ne\emptyset$ and
$\cod\ID'>r$.  Suppose $\ID'$ contains a subdiagram $\ID$ such that
$e\le3(r+1)$ and $s\ge r+1$.  Then the image of $Z(\ID)$ contains
$Y(\ID')$.  Hence $\cod Y(\ID')\ge s>r$.

Such a subdiagram $\ID$ exists by \cite{KP99, Lemma (4.4), p.~228}, if
$r\le7$, and we can easily extend the proof if $r=8$.  (A different
proof, valid for any $r$, is given in \cite{K--P}.)  We need only check
the case where $\ID'$ has only one root, say $R$ with weight $m'$.  If
$m'\ge5$, then take $\ID$ to consist of $R$ with weight 5 so that $e=15$
and $s=13$.  If $m'=4$, then $\ID'$ cannot have only one vertex since
$\cod\ID'>8$; hence, we may take $\ID$ to be ${\bf X}_{1,1}$ if $R$ is
followed by a vertex of weight 1, and to be ${\bf X}_{1,2}$ if $R$ is
followed by a vertex of weight 2 or more.  If $m'=3$, then $\ID'$ is
either $J_{l,j}$ or $E_{6l+j}$ where $l\ge2$; hence, we may take $\ID$
to be $J_{2,0}$ for which $e=12$ and $s=10$.  Finally, if $m'=2$, then
$\ID'$ is $A_k$ with $k>8$, and we may take $\ID$ to be $\IA_9$ so that
$e=15$ and $s=9$.

Third, let $\ID'$ be a diagram with $Y(\ID')\ne\emptyset$ and
$\cod\ID'\le r$.  Then $\deg\ID'\le3r$ by \cite{KP99, (4.3), p.~227}, since
$r\le8$ (a different proof, valid for any $r$, is given in \cite{K--P}).
Hence we may take $\ID'$ as $\ID$.  Then the image of $Z(\ID)$ contains
$Y(\ID')$.  Hence $\cod Y(\ID') \ge\cod\ID'$.  Thus Hypothesis (ii) of
Theorem (2.5) holds as well.

Therefore, by Theorem (2.5), the closure $\?{Y(r\IA_1)}$ is the support
of a nonnegative cycle $U(r)$, whose class $[U(r)]$ is given by a
certain expression.  We work it out in a moment.  (It is here alone
that we need the restriction $r\le8$.)  First, however, observe that
$U(r)$ is reduced; indeed, $Z(r\IA_1)$ is reduced (in fact, smooth) as
we proved above, so $U(r)$ is reduced by Lemma~(2.6).

Let $M$ be the subscheme of $Y$ parameterizing the curves that pass
through $g$ given general points; $M$ is an intersection of divisors,
one for each point, see the paragraph after (5.2.6) below.  Hence $M\cap
U(r)$ is reduced by the proof of \cite{KP99, Lemma (4.7), p.~232}, which
works virtually without change in the present setting.  So we have
  $$N_{g,r}=\int [M]\cdot [U(r)]$$
 once we've shown each point of $M\cap U(r)$ represents a curve $C$
that's irreducible.

Suppose some $C$ is reducible, say $C=C_1+C_2$.  Then each $C_i$ has
only ordinary nodes, say $r_i$ of them.  Moreover, the $C_i$ meet
transversally, say in $r_{12}$ points.  Then
        $$r=r_1+r_2+r_{12}.$$
 Let $d_i$ be the self-intersection number of $C_i$.  Then
        $$d=d_1+d_2+2r_{12}.$$

 Let $Y_i$ be the component of $C_i$ in the Hilbert scheme of $A$.  Then
every irreducible component of $Y_i(r_i\IA_1)$ is of dimension at least
$d_i/2-r_i+1$ (with equality if the appropriate $k_i$-ampleness holds).
Let $Y_i'$ be the component of $C_i$ in $Y_i(r_i\IA_1)$.  Summing
divisors induces a map $Y_1'\x Y_2'\to \?{Y(r\IA_1)}$.  Its fibers are
finite.  Hence
        $$\dim Y(r\IA_1)\ge (d_1/2-r_1+1)+(d_2/2-r_2+1).$$
  The right side is equal to $d/2-r+2$.  The left side is equal to
$d/2-r+1$.  Thus we have a contradiction.  So $C$ is irreducible.

Let's now work out the expression for $[U(r)]$.  First, note that
$w_i=0$ because $\Omega_A^1=0$ since $A$ is Abelian.  So each $b_q$
reduces to a certain polynomial in $v$.  So, to find $a_q:=\pi_*b_q$,
we must find $\pi_*v^a$ for $a\ge0$.

By definition, $v:=[D]$.  So (5.2.3) yields
        $$v = (1\x p)^*l+\pi^*h \where l:=c_1(\cP\ox p_1^*\cL)
        \and h:=c_1(\cO_Y(1)).$$
 Now, $\pi_*(1\x p)^*l=p^*{p_2}_*l$.  Hence the projection formula yields
        $$\pi_*v^a=\sum{a\choose i}(p^*{p_2}_*l^i)h^{a-i}.$$
 So we must compute ${p_2}_*l^i$ for $i\ge0$.

Let $\mu: A\x A\to A$ denote the group law.  Given $x\in A$, define $T_x\:
A\to A$ by $T_xy:=x+y$.  Finally, define $\phi\: A \to \hat A$ by
$\phi (x):=T_x^*\cL \ox \cL^{-1}$.  Then
  $$(1\x\phi)^*\cP = \mu^*\cL \ox p_1^*\cL^{-1}\ox p_2^*\cL^{-1}.$$
 according to \cite{M70, \p.151}.  Therefore,
        $$(1\x\phi)^*l=\mu^*c-p_2^*c \where c:=c_1(\cL).$$

Consider the Cartesian diagram
        $$\CD
         A\x A @>1\x\phi>> A\x\wh A \\
        @VV p_2 V            @VV p_2 V \\
        A          @>\phi>>    \wh A
  \endCD
  $$
 Note that $\phi^*{p_2}_*l^i = {p_2}_*(1\x\phi)^*l^i$.

The preceding two equations and the projection formula yield
  $$\phi^*{p_2}_*l^2 = {p_2}_*\mu^*c^2-
        2c{p_2}_*\mu^*c+c^2{p_2}_*[A\x A].$$
 Now, ${p_2}_*\mu^*c=0$ and ${p_2}_*[A\x A]=0$ by reason of dimension.
Also, $\int c^2=d$ by definition of $c$, $\cL$, and $d$.  Furthermore,
${p_2}_*\mu^*[x]=[A]$ for any $x\in A$, because
  $$\mu^{-1}x=\{\, (y,\,x-y)\mid y\in A\,\}.$$
 Hence $\phi^*{p_2}_*l^2 = d[A]$.  Since $[A]=\phi^*[\wh A]$,
therefore ${p_2}_*l^2 = d[\wh A]$.

 Similarly, we obtain
 $$\phi^*{p_2}_*l^3=-3dc \and \phi^*{p_2}_*l^4 =6dc^2.\eqno(5.2.5)$$

  Now, ${\phi}_*\phi^*z= (\deg\phi) z$ for any $z$ by
the projection formula, and $\deg \phi = d^2/4$ by \cite{M70,
\p.150}.  Taking $z:=l^3$ and $z:=l^4$, we therefore get
  $${p_2}_* l^3= -(12/d) {\phi}_*c\
        \and {p_2}_* l^4= (24/d){\phi}_*c^2.\eqno(5.2.6)$$
 Of course, ${p_2}_*l^i=0$ for $i\ne 2,3,4$ by reason of dimension.

We can now mechanically work out an expression for $[U(r)]$ as a linear
combination of $(p^*{\phi}_*c^i)h^{r-i}$ for $i=0,1,2$.  However,
 $N_{g,r}=\int [M]\cdot [U(r)]$.  So we must find $[M]$.  Well, given a
point $x\in A$, define $\iota_x\:Y\to F$ by $\iota_x(y):=(x,y)$.  Then
the divisor $\iota_x^{-1}D$ parameterizes the curves that pass through
$x$.  Owing to (5.2.3), the class $[\iota_x^{-1}D]$ is numerically the
same as $h$.  Hence $[M]$ is the same as $h^g$.  Therefore, owing to the
projection formula, we have to find $p_*h^{g+r-i}$ for $i=0,1,2$.

Recall that $Y:=\IP(\cQ^*)$.  So $p_*h^{g+r-i}=(-1)^is_{2-i}(\cQ)$ for
$i=0,1,2$ where the $s_j(\cQ)$ are the Segre classes.  Hence,
        $$p_*h^{g+r-2}=[\wh A]\and
          p_*h^{g+r-1}=-c_1(\cQ)\and
          p_*h^{g+r-1}=c_1(\cQ)^2-c_2(\cQ).$$
 It therefore remains to find $c_1(\cQ)$ and  $c_1(\cQ)^2$ and $c_2(\cQ)$.

Since $A$ is abelian, the Todd class of $p_2$ is trivial.  So the
Riemann--Roch theorem yields the following relation among the Chern
characters:
   $${\ch}({p_2}_*(\cP\ox p_1^*\cL))={p_2}_*({\ch}(\cP\ox p_1^*\cL)).$$
 Now, $\cQ:={p_2}_*(\cP\ox p_1^*\cL)$ by (5.2.1).  So, owing to (5.2.2),
the left side is equal to
        $$d/2+c_1(\cQ)+(c_1(\cQ)^2-2c_2(\cQ))/2.$$
 On the other hand, the right side is equal to ${p_2}_*(\sum l^i/i!)$.
Hence
        $${p_2}_*l^3=6c_1(\cQ)\and
        {p_2}_*l^4=12\bigl(c_1(\cQ)^2 - 2c_2(\cQ)\bigr).\eqno(5.2.7)$$
 So (5.2.6) yields
        $$c_1(\cQ)=-(2/d){\phi}_*c\and
        c_2(\cQ)=(1/2)c_1(\cQ)^2-(1/d){\phi}_*(c^2).\eqno(5.2.8)$$

To find $c_1(\cQ)^2$, use the formulas leading to (5.2.6).  Taking
$z:=c_1(\cQ)^2$ gives
 $$c_1(\cQ)^2=(4/d^2){\phi}_*\bigl({\phi}^*c_1(\cQ)\bigr)^2$$
 since $\phi^*z^2=(\phi^*z)^2$.  Now, ${\phi}^*c_1(\cQ)=(-d/2)c$ by
(5.2.7) and (5.2.5).  Hence
        $$c_1(\cQ)^2={\phi}_*(c^2).\eqno(5.2.9)$$

By definition, $\int c^2=d$.  So (5.2.9) and (5.2.8) yield
        $$\int c_1(\cQ)^2=d\and \int c_2(\cQ)=(d/2)-1. $$
 It is now a purely mechanical matter to derive the formulas in
Table~(5.1), and the proof is complete.\qed

 \lem3 Let $S$ be a smooth projective surface with numerically trivial
canonical bundle and with Picard number $1$.  Let $\cN$ be a line bundle
whose homology class is $m$ times the positive primitive class.  Set
$d:=\int c_1(\cN)^2$, and let $k\ge0$.  Assume either
 \rbit i that $m=1$ and $d>4(k+1)$ or
 \rbit ii that $m \ge2$ and $(m-1)d> m^2(k+1)$.

 \noindent Then $\cN$ is $k$-very ample.
 \pf
 If (ii) holds, then $4(m-1)d> 4m^2(k+1)$, and so, since $m^2\ge
4(m-1)$, then $d>4(k+1)$.  Hence, if either (i) or (ii) holds, then
$d>4(k+1)$.

Let $\cK$ be the canonical bundle.  Form $\cK^{-1}\ox\cN$, and to it,
apply Theorem 2.1 of Beltrametti--Sommese's \cite{BS88, p.~38}. Their
theorem implies that, since $d\ge4k+5$, either $\cN$ is $k$-very ample or
there exists an effective divisor $D$ such that
  $$\int c_1(\cN)[D]-(k+1)\le\int[D]^2\le\int c_1(\cN)[D]/2.\eqno(5.3.1)$$
 Suppose such a  $D$ exists, and let's derive a contradiction.

Say the class $[D]$ is $t$ times the positive primitive class.  Then the
second inequality in (5.3.1) becomes $t\le m/2$.  So $m\ge2$ since
$t\ge1$.  Hence Hypothesis (ii) applies, and so $(m-1)d> m^2(k+1)$.
However, the first inequality in (5.3.1) amounts to $t(m-t)d\le
m^2(k+1)$; whence, $(m-1)d\le m^2(k+1)$ because $m-1\le t(m-t) $ when
$1\le t\le m/2$.  Thus we have a contradiction, and the proof is
complete.\qed

 \rmk4 Similarly, we can enumerate those of the $C$ that lie in a given
linear equivalence class and pass through only $g-2$ general points.  In
fact, modified slightly, the proof of Theorem (5.2) yields the desired
enumeration, and shows it is valid when $r\le8$ and
        $$m>(3r+5)/2.\eqno(5.4.1)$$
 Conceivably, some $C$ are reducible, although the corresponding
dimension count shows that reducibility is not to be expected; compare
\cite{Gtt98, Rmk.~3.1, p.~528}.

 Indeed, the $C$ in the class are parameterized by a fiber of $Y/A$;
hence, the enumeration can be accomplished by computing the coefficient
of $\bigl[\wh A\,\bigr]$ in $p_*(h^{g-2}\cdot [U(r)])$.  Furthermore,
(5.4.1) implies $m\ge2$.  Also, (5.4.1) is equivalent to $2(m-1)>
3(r+1)$.  Now, $d/m^2$ is the self-intersection number of the primitive
class; so $d/m^2\ge2$.  Hence Lemma~(5.3) implies that, if $\cN$ is any
fiber of $\cP\ox p_1^*\cL$, then $\cN$ is $k$-very ample for
$k:=3(r+1)-1$.  The rest of the proof of validity is virtually the same.

Although Theorem (1.1) of \cite{KP99} yields the same formula on setting
$k$, $s$, and $x$ equal to 0, that theorem only asserts validity when
$m\ge3r$ and $\cO_A(C)=\cM^{\ox m}\ox\cN$ where $\cM$ is very ample and
$\cN$ is spanned.  On the other hand, that theorem does not require the
Picard number to be 1 nor the surface to be Abelian.  In any event, the
formula agrees with  G\"ottsche's Conjecture 2.4 in \cite{Gtt98, p.~526}.

Some condition like (5.4.1) is necessary.  Indeed, in \cite{D99,
Rmk.~2.3, p.~581}, Debarre considered the case in which $m$ is prime,
$d=2m^2$, and $g=2$, whence $r=m^2-1$.  In this case, G\"ottsche's
formula fails as G\"ottsche \cite{Gtt98, Rmk.~3.1, p.~528,} expected when
$Y(\infty)$ is nonempty.

On the other hand, the case $m=1$ is rather interesting and has
attracted some attention.  In this case, the method of Theorem (5.2)
shows that G\"ottsche's formula is valid if $A$ has Picard number 1 and
if $r\le8$ and $d>12(r+1)$, whereas Theorem (1.1) of \cite{KP99} asserts
nothing.  Using symplectic methods, Bryan and Leung \cite{BL98} showed
that the formula is valid for all $r$ and $d$ when $A$ is generic among
the Abelian surfaces for which the given homology class is algebraic.
Using complex analytic methods, in \cite{Gtt98, Thm.~3.2, p.~528},
G\"ottsche showed that the formula is valid for all $r$ when $g=2$ and
the homology class is a polarization of type $(1,n)$.  Independently and
somewhat differently, Debarre \cite{D99} proved the same result.
Earlier, see \cite{V95, Ex.~5.6, p.~521}, Schoen treated the case of a
polarization of type $(1,5)$ on a general Horrocks--Mumford Abelian
surface.

 \rmk5 In Lemma (5.3), if $S$ is a K3 surface, then in Condition (i) we may
replace $d>4(k+1)$ with $d\ge4k$.  Indeed, the proof is the same, except
that, instead of applying Theorem 2.1 of Beltrametti--Sommese
\cite{BS88, p.~38}, we apply Theorem 1.1 of Knutsen \cite{K01}, which
says that, since $S$ is a K3 surface and $d\ge4k$, either $\cN$ is
$k$-very ample or there exists an effective divisor $D$ such that
(5.3.1) holds.

  Similarly, if $S$ is an Enriques surface, then we may replace both (i)
and (ii) with the single condition that $d\ge 4(k+1)$.  Indeed, Theorem
1.2 of Knutsen's \cite{K01} asserts that, since $S$ is an Enriques
surface and $d\ge4(k+1)$, either $\cN$ is $k$-very ample or there exists
an nonzero effective divisor $D$ with nonpositive self-intersection; however,
the latter is impossible since $S$ has Picard number $1$ by hypothesis.

For these $S$, the method of proof of Theorem (5.2) shows that the
formula provided by Theorem (1.1) of \cite{KP99} is valid for more $m$.
Of course, some restriction on $m$ is necessary.  Indeed, in \cite{Ta82,
Ex. (3.13), p.~252}, Tannenbaum gave a simple example of a complete
linear system on a K3 surface $S$ such that $\cod Y(4\IA_1)<4$.

In Tannenbaum's example, $S$ is an arbitrary smooth quartic in $\IP^3$.
The system is the one cut by the quadrics; so it is parameterized by a
projective space $Y$ of dimension 9.  A general plane section of $S$ is
smooth.  So a general plane-pair section has two smooth components that
meet transversally in four points.  Hence $\dim Y(4\IA_1)\ge6$, and so
$\cod Y(4\IA_1)<4$.

Furthermore, if $S$ is generic, then its Picard group is generated by
$\cO_S(1)$ by the Noether--Lefschetz theorem.  In particular, the Picard
number is 1.  Also, if a quadric section is not reduced, then it must be
twice a plane section.  Since planes and quadrics are determined by
their sections, the quadric must be a double plane.  Hence $\dim
Y(\infty)=3$, and so $\cod Y(\infty)=6$.  Thus, unexpectedly, there are
infinitely many 4-nodal quadric sections through 5 general points, and
all are reduced.

  \serial{ajm}{Amer. J. Math.}
   \serial{cia}{Comm. in Alg.}
   \serial{cmp}{Comm. Math. Phys.}
   \serial{comp}{Comp. Math.}
   \serial{cnt}{Cont. Math.}
  \serial{duke}{Duke Math. J.}
   \serial{ens}{Ann. Sci. \'Ecole Norm. Sup.}
  \serial{imrn}{Int. Math. Res. Not.}
   \serial{invm}{Invent. Math.}
   \serial{indagmath}{Indag. Math.}
   \serial{ja}{J. Alg.}
   \serial{jag}{J. Alg. Geom.}
   \serial{jram}{J. Reine Angew. Math.}
   \serial{msm}{manuscripta math.}
   \serial{mathann}{Math. Ann.}
   \serial{mnach}{Math. Nachr.}
   \serial{mz}{Math. Z.}
   \serial{pjm}{Pac. J. Math.}
   \serial{pams}{Proc. Amer. Math. Soc.}
   \serial{plms}{Proc. Lond. Math. Soc.}
   \serial{rspad}{Rend. Sem. Mat. Univ. Padova}
   \serial{sms}{Selecta Math. Sovietica}
   \serial{tams}{Trans. Amer. Math. Soc.}

\references

BS88
      M. Beltrametti and A. J. Sommese,
      Zero cycles and $k$-th order embeddings of smooth projective
      surfaces,
      In Problems in the theory of surfaces and their classification
      (Cortona, 1988), Sympos. Math., XXXII, Academic Press, London,
      1991, 33--48.

B--L99
     C. Birkenhake and H. Lange,
     ``Complex tori,"
     Birkh\"auser, Progress in mathematics {\bf 177}, 1999.

BL98
     J. Bryan and N. C. Leung,
     Generating functions for the number of
     curves on abelian surfaces,
      \duke 99 1999 311--28.

CH98
     L. Caporaso and J. Harris,
     Counting plane curves of any genus,
     \invm 131 1998 345--92.

C66
    A. Cayley, On the theory of involution,
    {\it Trans. Cambridge Phil. Soc.,} {\bf XI. Part I} (1866), 21--38.
    = ``Coll. Math. Papers of A. Cayley,'' {\bf  V},
    [348], 1892, pp.~295--312.

Ch97
    Y. Choi,
    On the degree of Severi varieties,
    preprint 1997 (available from http://math.ucr.edu
    /$\tilde{}\,$ychoi/\break paper.html).

C86
H. Clemens,
  Curves on higher-dimensional complex projective manifolds,
  in ``Proc. International Cong. Math., Berkeley, 1986,'' pp.~634--40

Co70
  L. Comtet,
  ``Analyse combinatoire I,''
  SUP {\bf 4}, Presses Universitaires de France, 1970.

CK99
     D. Cox and S. Katz,
     ``Mirror symmetry and algebraic geometry,''
     Math. Surveys and Monographs, Vol. {\bf 68}, AMS 1999.

D99
     O. Debarre,
     On the Euler characteristic of generalized Kummer varieties,
     \ajm 121 1999 577--86.

Gtt98
      L. G\"ottsche,
      A conjectural generating function for numbers of
      curves on surfaces,
      \cmp {\bf 196} 1998 523--33.

GL96
      G.-M. Greuel and C. Lossen,
      Equianalytic and equisingular families of curves on surfaces,
      \msm 91 1996 323--42.

GLS97
 G.-M. Greuel, C. Lossen, and E. Shustin,
 New asymptotics in the geometry of equisingular families of curves
  \imrn 13 1997 595--611.

HP95
      J. Harris and R. Pandharipande,
      Severi degrees in cogenus 3,
      alg-geom/9504003.

JK96
      T. Johnsen and S. Kleiman,
      Rational curves of degree at most $9$ on a general quintic
      threefold,
      \cia {\bf 24}(8) 1996 2721--53.

TdJ00
   T. de Jong,
   Equisingular deformations of plane curve and of sandwiched
singularities,
   arXiv: math.AG/\break0011097.

Ka86
  S. Katz,
  On the finiteness of rational curves on quintic threefolds,
  \comp 60 1986 151--62.

Ka92
  S. Katz,
  Rational Curves on Calabi-Yau Threefolds,
  in ``Essays on mirror manifolds,'' S.-T. Yau (ed.), Int. series in
math. physics, International Press (1992), pp. 168--80.

Kl74
   S. L. Kleiman,
  The transversality of a general translate,
  \comp 28 1974 287--97.

KP99
     S. Kleiman and R. Piene,
     Enumerating singular curves on surfaces,
     in ``Algebraic geometry --- Hirzebruch 70,'' \cnt 241 1999 209--38
     \ (corrections and revision in math.AG/9903192).

K--P
     S. Kleiman and R. Piene,
     Node polynomials for curves on surfaces,
     to appear.

K01
     A. L. Knutsen,
     On $k$th order embeddings of K3 surfaces and Enriques surfaces,
     \msm 104 2001 211--37.

Los98
   C. Lossen,
   ``The geometry of equisingular and equinalytical families of curves on
   a surface,''
   Dr. dissertation, Universit\"at Kaiserslautern, 1998.

Mat91
   J. F. Mattei,
   Modules de feuilletages holomorphes singuliers: I
   \'equisingularit\'e,
   \invm 103 1991 297--325.

M70
     D. Mumford,
     ``Abelian Varieties," Oxford University Press, 1970.

NV97
 A. Nobile and O. E. Villamayor,
 Equisingular stratifications associated to families of planar ideal,
 \ja 193 1997 239--59.

Ran89
   Z. Ran,
   Enumerative geometry of singular plane curves,
   \invm 97 1989 447--65.

R67
     S. Roberts,
	Sur l'ordre des conditions pour la coexistence des \'equations
	alg\'ebriques \`a plusieurs variables,
     \jram 67 1867 266--78.

R75
     S. Roberts,
 On a simplified method of obtaining the order of algebraical conditions,
  \plms {} 1875 101--13.

S65
    G. Salmon,
    ``A treatise on the analytic geometry of three dimensions,''
    2nd edition,
    Hodges, Smith, and Co., Dublin 1865

S79
 G. Salmon,
	``Higher plane curves," 3rd edition 1879, Chelsea reprint.

St48
 J. Steiner,
  Allgemeine Eigenschaften der algebraischen Curven,
  {\it Berlin. Ber.} 1848, 310--15 $=$ \jram 47 1853 1--6\ $=$
``Ges. Werke,'' herausg. von K. Weierstrass, Berlin 1882, {\bf 2}, 495--500.

Ta82
     A. Tannenbaum,
     Families of curves with nodes on K3-surfaces,
     \mathann 260 1980 239--53.

Te73
 B. Teissier,
 Cycles \'evanescents, sections planes et conditions de Whitney,
 in ``Singularit\'es \`a Carg\`ese," Ast\'erisque {\bf 7--8} (1973),
285--362.

V95
     I. Vainsencher,
     Enumeration of $n$-fold tangent hyperplanes to a surface,
     \jag {\bf 4} 1995 503--26.

Wahl74
  J. M. Wahl,
 Equisingular deformations of plane algebroid curves, \tams
 193 1974 143--70.

Wall84
   C. T. C. Wall,
   Notes on the classification of singularities,
   \plms 48 1984 461--513.

Z82
    O. Zariski,
    Dimension-theoretic characterization of maximal irreducible
    algebraic systems of plane nodal curves of a given order $n$ and
    with a given number $d$ of nodes,
    \ajm 104 1982 209--26.

\endreferences

\medskip
 \centerline{\vrule width 3 true cm height 0.4pt depth 0pt}

 \bigskip

\eightpoint
%  \centerline{Mathematics Department, Room {\sl 2-278} MIT,
  \centerline{Department of Mathematics, Room {\sl 2-278} MIT,
   {\sl77} Mass Ave, Cambridge, MA {\sl02139-4307}, USA}
 \smallskip
  \centerline{E-mail: {\tt kleiman\@math.mit.edu}}
 \medskip
  \centerline{Department of Mathematics, University of Oslo,
   PO Box {\sl1053}, Blindern,  NO-{\sl0316} Oslo, Norway}
    \smallskip
  \centerline{E-mail: {\tt ragnip\@math.uio.no}}
\bye